\documentclass[10pt]{article}
\usepackage{amsfonts,epsfig,a4wide}

\newtheorem{Th}{Theorem}

\newtheorem{Lemma}{Lemma}
\newtheorem{Cor}{Corollary}
\newtheorem{Def}{Definition}
\newtheorem{Rem}{Remark}

\newcommand{\eqdef}{\stackrel{{\rm def}}{=}}

\newcommand{\trace}{\mbox{\rm trace}}

\newcommand{\diag}{\mbox{\rm diag}}

\newcommand{\Id}{\mbox{\rm Id}}

\newcommand{\weg}[1]{}

\begin{document}
\title{Projective Lichnerowicz-Obata Conjecture}
\author{Vladimir S. Matveev\thanks{Mathematisches Institut, Universit\"at
Freiburg, 79104 Germany \
 matveev@email.mathematik.uni-freiburg.de}}
\date{}
\maketitle

\begin{abstract}
We solve two classical conjectures by showing that
 if  an action of a connected  Lie group  on a complete
 Riemannian manifold preserves   the geodesics
(considered as unparameterized curves), then   the metric has
constant positive sectional curvature, or the group acts by affine transformations.
\end{abstract}

\section{Introduction} \label{introduction}
\subsection{Results}

\begin{Def}\hspace{-2mm}{\bf .} Let $(M_i^n, g_i)$, $i=1,2$,  be   smooth Riemannian manifolds.

A diffeomorphism $F:M_1^n\to M_2^n$ is called {\bf projective}, if
it takes unparameterized geodesic of $g_1$  to geodesics of $g_2$.
A projective diffeomorphism of a Riemannian manifold is called
a {\bf projective transformation}.

A diffeomorphism $F:M_1^n\to M_2^n$ is called {\bf affine}, if it
takes the Levi-Civita connection of   $g_1$  to  the Levi-Civita
connection of $g_2$. An affine  diffeomorphism of a Riemannian
manifold is called an {\bf affine  transformation.}
\end{Def}

\begin{Th}[Projective Lichnerowicz  Conjecture]\hspace{-2mm}{\bf .} \label{main}
Let a connected  Lie group $G$  act on a complete
 connected  Riemannian manifold $(M^n, g)$ of dimension
 $n\ge 2$   by projective
transformations.  Then, it  acts by affine transformations, or $g$
has  constant positive sectional  curvature.
\end{Th}

The Lie groups of affine transformations of
complete Riemannian manifolds are well understood,
  see
 \cite{lichnerowicz}: suppose a  connected Lie group acts on
 a  simply-connected complete
 $(M^n,g)$   by affine transformations.
Then, there exists a Riemannian decomposition $(M^n,
g)=(M_1^{n_1}, g_1)+ (\mathbb{R}^{n_2}, g_{{euclidean}})$ of the
manifold into the direct sum of a Riemannian manifold $(M_1^{n_1},
g_1)$ and  Euclidean space $(\mathbb{R}^{n_2}, g_{{euclidean}})$
such that
 the group acts   componentwise; it acts by isometries  on
$(M_1^{n_1}, g_1)$ and by compositions
 of   linear transformations  and   parallel translations
on $\mathbb{R}^{n_2}$.  In particular, every   connected Lie group
of affine transformations of  a closed manifold  consists of
isometries \cite{yano1}.
Thus,  the following statement is a direct consequence of
Theorem~\ref{main}

\begin{Cor}[Projective Obata Conjecture]\hspace{-2mm}{\bf .} \label{obata}
Let a connected  Lie group $G$  act on a closed connected
Riemannian manifold $(M^n, g)$ of dimension
 $n\ge 2$   by projective
transformations.  Then, it acts by isometries, or $g$ has constant
positive  sectional curvature.
\end{Cor}

\voffset=-2.2 truecm  \textheight=9.8 truein 

\evensidemargin -3 in

Any   connected simply-connected Riemannian manifold of constant
positive  sectional  curvature is  a  round sphere. All projective
transformations of the round sphere are  known (essentially, since
Beltrami \cite{Beltrami});
 so that Theorem~\ref{main}  closes the theory of nonisometric
infinitesimal  projective transformations of complete  manifolds.

\subsection{History}

 The theory of  projective transformations
has a long and  fascinating history.
First  nontrivial examples of projective transformations
 are due to  Beltrami \cite{Beltrami}. We
describe their natural  multi-dimensional generalization.  Consider
the sphere $$ S^n\eqdef \{(x_1,x_2, ...,x_{n+1})\in R^{n+1}: \ x_1^2+x_2^2+...+x_{n+1}^2=1\}
$$ with the restriction of the Euclidean metric   and
 the mapping $a:S^n\to S^n$  given by
$a:v\mapsto \frac{A(v)}{\|A(v)\|}$, where $A$ is an
arbitrary
non-degenerate linear  transformation of $R^{n+1}$.

The mapping is clearly a diffeomorphism taking geodesics to geodesics.
Indeed, the geodesics of  $g$ are great circles
(the intersections of planes that go through the origin with the
 sphere). Since $A$ is linear, it
 takes
planes to  planes. Since the normalization $w\mapsto
\frac{w}{\|w\|}$ takes punctured   planes to their intersections with the
sphere,  $a$ takes  great circles to great circles. Thus,  $a$ is a
projective transformation.
 Evidently, if $A$ is not proportional to an
  orthogonal transformation,  $a$ is not affine.

Beltrami investigated separate examples of projective transformations. One of    the first important papers  on smooth families of projective
transformations is due to Lie, see  \cite{Lie}.  Lie formulated the
problem of finding metrics (on surfaces) whose groups of projective
transformations are  bigger than the groups of isometries  ("Lie
Problem" according to Fubini), and solved it assuming that the
groups are  big enough. In the case  when the manifold is complete, this problem was
formulated in Schouten \cite{Schouten}.

The local theory of projective transformations was well understood
thanks to efforts of several mathematicians, among them   Dini \cite{Dini}, Levi-Civita \cite{Levi-Civita},
Fubini \cite{Fubini1}, Eisenhart \cite{Eisenhart}, Cartan
\cite{cartan},  Weyl \cite{Weyl2} and Solodovnikov \cite{solodovnikov1}.
We will recall  their
results in Theorems~\ref{LC},\ref{maximal},\ref{solodovnikov},\ref{fubini}.

 A basic
philosophical  idea  behind these results
can be described as follows (see, for example, \cite{Weyl1}): the Universe can
be
 explained by its  infinitesimal structure, and    this infinitesimal structure is invariant
 with respect to a group of transformations.

Weyl studied projective transformations  on the tensor level and
found a number of  tensor reformulations.
 He constructed
 the so-called  projective Weyl tensor  $W$ \cite{Weyl2}
which is invariant with respect to projective transformations.
We will recall the  definition of $W$  in Section~\ref{Vk}, and will
 use it is Section~\ref{3.5}.
E. Cartan \cite{cartan},  T. Y. Thomas \cite{thomas},
J. Douglas \cite{douglas}  and  A. Lichnerowicz et al \cite{lichnerowicz2}
  studied  groups of
 projective  transformations  on the
level of   affine connections,   sprays and natural Hamiltonian  systems. They introduced the
  so-called projective  connection and  Thomas  projective
parameters, which are invariant with respect to
  projective transformations.

Theorem~\ref{main} and Corollary~\ref{obata} are known in mathematical folklore as the
Lichnerowicz  and Obata conjectures, respectively, although Lichnerowicz and Obata never formulated them explicitly.  In  several papers
(see,  for example, \cite{nagano,Yamauchi1,hasegawa}), they were formulated as    "well  known classical conjectures".

Maybe  the name "Lichnerowicz-Obata conjecture" appeared because
of the  similarity with the conformal Lichnerowicz conjecture
(proved by Obata \cite{Obata}, Alekseevskii \cite{alekseevskii}
and Ferrand~\cite{ferrand}).

Note  that, in the time of Lichnerowicz and Obata, projective and
conformal transformations were studied by the same people by the
same methods (see, for example, \cite{Couty1,YN}):
the tensor equations for conformal and projective
 infinitesimal transformations are very similar.

Projective transformations  were  extremely popular objects of
investigation in 50th--80th. One of the reason for
it is their possible applications in physics, see, for example,
\cite{petrov,davis}. One may consult the
surveys~\cite{mikes} (more geometric one)  and \cite{aminova}
(from the viewpoint of physics), which contain  more than 500
references.

Most results on projective transformations require additional
geometric assumptions written as a tensor equation. For example,
Corollary~\ref{obata} was proved under the assumption that
 the metric is
Einstein \cite{Couty},   K\"ahler \cite{Couty}, Ricci-flat \cite{nagano2},  or has constant
scalar curvature \cite{Yamauchi1}.

 An  important result which does
not require additional tensor assumptions is due to Solodovnikov.
He proved the  Lichnerowicz
 conjecture under the assumptions
that the dimension of the manifold is greater than two and    that all  objects
(the metric, the manifold, the projective transformations) are
real-analytic. The statement itself is in \cite{solodovnikov3},
but the technique was mostly developed in \cite{solodovnikov1,solodovnikov2,solodovnikov2.5,solodovnikov2.6}.
In Section~\ref{Vk},  we
will review the results of Solodovnikov  we need
for proving Theorem~\ref{main}. We will also use certain results of Solodovnikov in Section~\ref{3.5}.

Both assumptions are important for Solodovnikov's methods.
Solodovnikov's technique  is  based on a very accurate analysis of
the behavior of the curvature tensor under projective
transformations and completely fails in dimension two,
see Theorem~\ref{solodovnikov}  and Examples 1,2 from Section~\ref{dim2}. Also,
 real analyticity is extremely  important  for his methods: all
 his global statements are based on it.

The new techniques  which allow us to prove the Lichnerowicz-Obata
conjecture were introduced in
\cite{MT,TM,dedicata,quantum,physica}: the main observation is
that the existence of  projective diffeomorphisms   allows one to
construct commuting integrals for the geodesic flow, see
 Theorem~\ref{integrability} in Section~\ref{2.2}. This observation  has been
used quite successfully in finding   topological obstruction that
prevent a closed manifold from possessing non-isometric projective
diffeomorphisms
 \cite{ERA, short,starrheit,topology,hyperbolic}.

\subsection{ Counterexamples to Theorem~\ref{main}, if one of the assumptions is
omitted }\label{contrexamples}  All assumptions in
Theorem~\ref{main} are important:

If the Lie group is not connected, a counterexample to
Theorem~\ref{main} is possible only if the group is discrete. In
this case, a counterexample exists  already in dimension two:
consider the torus $T^2:={\mathbb{R}^2}/_{\mathbb{Z}^2}$ with the
standard coordinates $x,y\in (\mathbb{R} \  \textrm{mod}  \ 1)$
and a positive smooth nonconstant function $f: (\mathbb{R} \
\textrm{mod}\ 1 )\ \to \mathbb{R}$ such that the metric
\begin{eqnarray*}
ds^2&:=&\left(f(x)-\frac{1}{f(y)}\right)\left(\sqrt{f(x)}\, dx^2+
\frac{1}{\sqrt{f(y)}}\, dy^2\right)
\end{eqnarray*}
is positive definite. Then, the diffeomorphism $F:T^2\to T^2$
given by $F(x,y):=(y,x)$ takes the metric   to
\begin{eqnarray*}
&&\left(f(y)-\frac{1}{f(x)}\right)\left(\frac{\sqrt{f(x)}}{f(x)}\, dx^2 +
\frac{f(y)}{\sqrt{f(y)}}\, dy^2
\right),
\end{eqnarray*}
and
 is a projective
transformation by Levi-Civita's Theorem~\ref{LC}.

If the manifold  is not complete, the first counterexamples can be
found in \cite{Lie}. In Section~\ref{2proof}, at the end of the
proof of Theorem~\ref{2dim},   we will essentially construct
generalizations of Lie's examples for every dimension $\ge 2$.

If the manifold is not connected, one can construct noninteresting
counterexamples as follows: the manifold has two connected
components. The first component is    the round sphere, where the
group $GL$ acts by projective transformations as in Beltrami's
example. The other component is something different from the
sphere where the group $GL$  acts identically.

If the manifold is 1-dimensional, every diffeomorphism is a
projective transformation and only homotheties are affine
transformations.

\subsection{Acknowledgements}

I would like to thank D.   Alekseevskii for formulation of the
problem,  V.  Bangert,   A. Bolsinov,  A. Fomenko,  I. Hasegawa,
M. Igarashi,   K. Kiyohara,    O. Kowalsky and K. Voss for useful
discussions,  DFG-programm 1154 (Global Differential Geometry) and
Ministerium f\"ur Wissenschaft, Forschung und Kunst
Baden-W\"urttemberg  (Elitef\"orderprogramm Postdocs 2003) for
partial financial support.

\section{ Preliminaries:  BM-structures,
    integrability,
 and Solodov\-nikov's $V(K)$ spaces
 } \label{preliminaries}

The  goal of this  section is to formulate  classical  and new
tools for  the proof of Theorem~\ref{main}. In
Sections~\ref{BM},\ref{integrals},
 we introduce the notion
"BM-structure" and explain its  relations to  projective
transformations and integrability; these are  new instruments of
the proof.  In Section~\ref{Vk},  we formulate in a convenient form
 classical results of Beltrami, Weyl,
Levi-Civita, Fubini, de Vries and Solodovnikov.    We will
actively use these results  in Sections~\ref{global},\ref{proof}.

\subsection{BM-structure} \label{BM}
Let $(M^n, g)$ be a Riemannian manifold of dimension $n\ge 2$.
\begin{Def}\hspace{-2mm}{\bf .}  \label{bm} A {\bf BM-structure} on   $(M^n, g)$ is a smooth
self-adjoint $(1,1)$-tensor  $L$  such
that,  for every point $x\in M^n$, for every vectors $u,v,w\in
T_xM^n$, the following equation holds:
\begin{equation}\label{condition}
g ((\nabla_u L) v, w)=
\frac{1}{2} g (v,u)\cdot d\trace_L(w)+
\frac{1}{2} g (w,u)\cdot d\trace_L(v),
\end{equation}
where $\trace_L$ is the trace of $L$.
\end{Def}

The set of all BM-structures on  $(M^n,g)$ will be denoted by
${\cal B}(M^n,g)$. It is a linear vector space. Its dimension is at least
one,  since  the identity tensor $\Id\eqdef
\diag(1,1,1,...,1)$ is always a BM-structure.

\begin{Def}\hspace{-2mm}{\bf .} Let $g, \ \bar g$ be Riemannian
metrics on $M^n$. They are {\bf projectively equivalent}, if they
have the same (unparameterized) geodesics.
\end{Def}

The relation between BM-structures and projectively equivalent metrics
is given by

\begin{Th}[\cite{benenti}]\hspace{-2mm}{\bf .} \label{th1}
 Let $g$ be a Riemannian metric. Suppose  $L$ is a
self-adjoint positive-definite $(1,1)$-tensor. Consider the
 metric $\bar g$ defined by
\begin{equation}\label{bg}
\bar g(\xi,\eta)=\frac{1}{\det(L)}
g(L^{-1}(\xi), \eta)
\end{equation}
(for every tangent vectors $\xi$ and $\eta$ with the common foot
point.)

Then,
the metrics $g$ and $\bar g$ are projectively equivalent,
 if and only if  $L$ is a  BM-structure on $(M^n,g)$.
\end{Th}

Equivalent form of this theorem  is

\begin{Cor}[\cite{benenti}]\hspace{-2mm}{\bf .} \label{th2}
 Let $g$,  $\bar g$ be  Riemannian metrics on $M^n$.

Then, they are projectively equivalent,  if and only if the tensor
 $L$ defined by
\begin{eqnarray}
L^i_j &\eqdef & \left(\frac{\det(\bar g)}{\det(g)} \right)^\frac{1}{n+1}
\sum_{\alpha=1}^n \bar g^{i\alpha} g_{\alpha j} \label{l}
\end{eqnarray}
is  a  BM-structure on $(M^n,g)$.
\end{Cor}

A one-parametric group of projective transformations of $(M^n, g)$
gives us a one parametric family  of BM-structures. Its derivative
is also a BM-structure:

 \begin{Th}[Infinitesimal version of Theorem~\ref{th1}]\hspace{-2mm}{\bf .} \label{projective}
  Let $F_t$, $t\in \mathbb{R}$,    be a
 smooth one-parametric family of projective transformations of  $(M^n, g)$.
 Consider  the (1,1)-tensor $A$  given by
\begin{equation}\label{A}
A^i_j\eqdef \sum_{\alpha=1}^ng^{i\alpha}\left({\cal L} g\right)_{\alpha j},
\end{equation}
 where ${\cal L} g$ denotes the Lie derivative with
respect to $F_t$ (so that $\left({\cal L} g\right)_{\alpha
j}=-\frac{d}{dt}\left((F_t^*g)_{\alpha j}\right)_{|t=0}$), and
$g^{i j}$ is  the inverse tensor to $g_{i j}$. Then,
$A-\frac{1}{n+1}\trace_A\cdot \Id$ is a BM-structure on $(M^n
,g)$.
 \end{Th}

\noindent {\bf Proof:} For every $t$, let us denote by $g_t$ the pull-back
$F^*_tg$. Fix a point $x\in M^n$  and a coordinate system at $T_xM^n$. Then,
we can think that $g_t$ and $g$ are matrixes. Clearly, $g_0=g$.
Since $F_t$ consists of projective transformations, by Corollary~\ref{th2}, for every $t\in \mathbb{R}$, the tensor
$$
L_t\eqdef \left(\frac{\det(g_t)}{\det(g)}\right)^{\frac{1}{n+1}}g^{-1}_tg
$$
satisfies the equation
\begin{equation}
g ((\nabla_u L_t) v, w)=
\frac{1}{2} g (v,u)\cdot d\trace_{L_t}(w)+
\frac{1}{2} g (w,u)\cdot d\trace_{L_t}(v)
\end{equation}
for every $u,v,w\in T_xM^n$.
Differentiating this equation by $t$ and substituting $t=0$, we obtain

$$
g \left(\left(\nabla_u \left(\frac{d}{dt}L_t\right)_{|t=0}\right) v, w\right)
=
\frac{1}{2} g (v,u)\cdot d\trace_{\left(\frac{d}{dt}L_t\right)_{|t=0}}(w)+
\frac{1}{2} g (w,u)\cdot d\trace_{\left(\frac{d}{dt}L_t\right)_{|t=0}}(v),
$$
so that $\left(\frac{d}{dt}L_t\right)_{|t=0}$ is a BM-structure on $M^n$.
Now, let us calculate it.
\begin{eqnarray*}
\left(\frac{d}{dt}L_t\right)_{|t=0} &=&
\left(\frac{\det(g_t)}{\det(g)}\right)_{|t=0}^{\frac{1}{n+1}}\left(\frac{d}{dt}g^{-1}_t\right)_{|t=0}g+\left(\frac{d}{dt}\left(\frac{\det(g_t)}{\det(g)}\right)^{\frac{1}{n+1}}\right)_{|t=0}\left(g_t^{-1}g\right)_{|t=0}\\
&= & \left(\frac{d}{dt}g^{-1}_t\right)_{|t=0}g+\left(\frac{d}{dt}\left(\frac{\det(g_t)}{\det(g)}\right)^{\frac{1}{n+1}}\right)_{|t=0}\Id\\
&=&-\left(g_t^{-1}\right)_{|t=0} \left(\frac{d}{dt}g_t\right)_{|t=0}
\left(g_t^{-1}\right)_{|t=0}g +\frac{1}{n+1}\left(\frac{d}{dt}
\frac{\det(g_t)}{\det(g)}\right)_{|t=0}
 \weg{\left(\frac{\det(g_t)}{\det(g)}\right)_{|t=0}} \Id\\
&=& g^{-1}{\cal L}g -\frac{1}{n+1}\trace_{g^{-1}{\cal L}g}
 \Id.
\end{eqnarray*}

Thus, $L:=A-\frac{1}{n+1}\trace_A\Id$ is a BM-structure. Theorem~\ref{projective} is proved.

For use in future we recall  one more property of BM-structures:

\begin{Th}[\cite{benenti},\cite{japan}]\hspace{-2mm}{\bf .} \label{nijenhuis}
The Nijenhuis torsion of a BM-structure vanishes.
\end{Th}

\subsection{Integrals for geodesic flows of metrics admitting
BM-structure} \label{integrals}\label{2.2}

Objects similar to    BM-structures on  Riemannian manifolds
appear  quite  often in the theory of integrable systems (see, for
example \cite{Benenti1,Benenti2,Benenti3,I,Crampin}). The relation between
BM-structures and    integrable geodesic flows is observed on the
level of projective equivalence in \cite{MT}, on the level of
projective transformations in \cite{Topalov} and is as  follows:

Let $L$ be a  self-adjoint  \ $(1,1)$-tensor on $(M^n,g)$.
Consider the family $S_t$, $t\in \mathbb{R}$, of $(1,1)$-tensors
\begin{equation}\label{st}
 S_t\eqdef \det(L - t\ \mbox{\rm Id})\left(L-t\ \mbox{\rm Id}\right)^{-1}.
 \end{equation}
\begin{Rem}\hspace{-2mm}{\bf .}
Although $\left(L-t\ \Id\right)^{-1}$ is not defined for
$t$
lying
in the spectrum of $L$, the tensor  $S_t$  is well-defined
for every  $t$.  Moreover,   $S_t$
is a  polynomial
 in $t$ of degree $n-1$
with coefficients being  (1,1)-tensors.
\end{Rem}
\noindent We will identify the tangent and  cotangent bundles of $M^n$ by $g$.
This identification allows us to transfer the
 natural  Poisson structure from  $T^*M^n$ to   $TM^n$.

\begin{Th}[\cite{TM, MT, dedicata,Topalov}]\hspace{-2mm}{\bf .}
\label{integrability}
 If $L$ is a BM-structure,
then, for every  $t_1,t_2\in \mathbb{R}$, the functions
\begin{equation}\label{integral}
I_{t_i}:TM^n\to \mathbb{R}, \ \ I_{t_i}(v)\eqdef g(S_{t_i}(v),v)
\end{equation}
are commuting integrals for the geodesic flow
 of  $g$.
\end{Th}

 Since $L$ is  self-adjoint,  its   eigenvalues are real. At every point $x\in M^n$,
let
 us denote by  $\lambda_1(x)\le ... \le \lambda_n(x)$  the eigenvalues
 of $L$ at the point.

\begin{Cor}\hspace{-2mm}{\bf .}\label{ordered12}
Let $(M^n, g)$ be a   connected Riemannian manifold such that
every two points can be connected by a geodesic. Suppose $L$ is a
BM-structure on $(M^n, g)$.
 Then,   for every  $i\in \{1,... ,n-1\} $, for every  $x,y\in M^n$, the
 following statements  hold:
\begin{enumerate}
\item $\lambda_i(x)\le \lambda_{i+1}(y)$.

\item  If $\lambda_i(x)< \lambda_{i+1}(x)$,
then $\lambda_i(z)< \lambda_{i+1}(z)$
for almost every point $z\in M^n$.


\end{enumerate}
\end{Cor}

A slightly different version of this corollary was proved in
\cite{hyperbolic,dedicata}. We will need the technique of the
proof later, and,  therefore,  repeat the proof in
Section~\ref{proofordered}.

At every point $x\in M^n$, denote by $N_L(x)$ the number of
different eigenvalues  of the BM-structure $L$ at $x$.

\begin{Def}\hspace{-2mm}{\bf .}
A point $x\in M^n$  will be called {\bf typical} with respect to the
BM-structure $L$, if $$ N_L(x)=\max_{y\in M^n}N_L(y). $$
\end{Def}

\begin{Cor}\hspace{-2mm}{\bf .}\label{typical}
Let  $L$ be a BM-structure on  a   connected Riemannian manifold
$(M^n, g)$.
 Then,  almost every point of $M^n $ is typical with respect to $L$.
\end{Cor}

\noindent {\bf Proof: } Consider points $x,y\in M^n$ such that $x$
is typical. Our goal is to prove that almost every point in a
small ball around $y$ is typical as well. Consider a path
$\gamma\in M^n$ connecting $x$ and $y$. For every point $z\in
\gamma$, there exists $\epsilon_z>0$ such that the open ball with
 center in $z$ and radius $\epsilon_z$ is convex. Since $\gamma$
is compact, the union  of a finite number of such balls contains
the whole path $\gamma$. Therefore, there exists a finite sequence
of  convex balls $B_1,B_2,...,B_m$ such that
\begin{itemize}
\item $B_1$ contains $x$.
\item $B_m$ contains $y$.
\item For every $i=1,...,m-1$, the intersection  $B_i\bigcap
B_{i+1}$ is not empty.
\end{itemize}
Since the balls are convex, every two points of every ball can be
connected by a geodesic. Using that, for a fixed $i$,  the set
$\{x\in M^n:\ \lambda_i(x)<\lambda_{i+1}(x)\}$
 is evidently open, by Corollary~\ref{ordered12},
almost every point of $B_1$ is typical. Then, there exists a
typical point in  the ball $B_2$. Hence, almost all points of
$B_2$ are  typical. Applying this argumentation $m-2$  times, we
obtain that almost all points of $B_m$ are  typical. Corollary is
proved.

\subsection{Projective Weyl tensor, Beltrami Theorem,
Levi-Civita's Theorem and Solodovnikov's $V(K)$ metrics} \label{Vk}
Let $g$ be a Riemannian metric on $M^n$ of dimension $n$.
Let $R_{jkl}^i$ and $R_{ij}$ be the curvature and the Ricci tensors of $g$.
The tensor
\begin{equation}\label{pw}
W_{jkl}^i:=R_{jkl}^i-\frac{1}{n-1}\left(\delta_l^i \, R_{jk}- \delta_k^i\, R_{jl}\right)
\end{equation}
is called the {\it projective Weyl} tensor.

\begin{Th}[\cite{Weyl2}]\hspace{-2mm}{\bf .} \label{weyl}
 Consider   Riemannian  metrics $g$ and $\bar g$ on $M^n$.
  Then, the
 following statements hold:

\begin{enumerate}
\item If the metrics are   projectively equivalent, then their projective Weyl tensors coincide.
\item Assume  $n\ge 3$. The projective Weyl tensor of $g$  vanishes if and only if
the sectional curvature  of $g$ is constant.
\end{enumerate}
\end{Th}

\begin{Cor}[Beltrami Theorem   \cite{Eisenhart}]\hspace{-2mm}{\bf
.} \label{belt} Projective diffeomorphisms take metrics of constant sectional
curvature to   metrics of constant  sectional curvature.
\end{Cor}

Formally speaking, Corollary~\ref{belt} follows from Theorem~\ref{weyl} for
dimensions greater than two only. For dimension two, Corollary~\ref{belt} was
proved by Beltrami himself in \cite{Beltrami}.

\begin{Cor}\hspace{-2mm}{\bf .} \label{bel} Let $F:M_1^n\to M_2^n$ be
a  projective diffeomorphism between complete  Riemannian manifolds
 $(M_i^n,g_i)$, $i=1,2$, of dimension $n\ge 2$.
Then, if $g_1$ has constant negative sectional
curvature, $F$ is a homothety.  If $g_1$ is flat,  $F$ is affine.
\end{Cor}

Corollary~\ref{bel}  is a mathematical folklore. Unfortunately, we did not
find a classical reference for  it.  If $g_1$ is flat,
Corollary~\ref{bel} can be found  in every good textbook
 on  linear algebra. If the
curvature of $g_1$ is negative,
 under the assumption that the  dimension  is two,
Corollary~\ref{bel} was proved in \cite{bonahon}. Case  arbitrary dimension
trivially follows from dimension two, since in every two-dimensional
direction there exists  a totally geodesic complete submanifold.

 In view of Theorem~\ref{th1}, the next   theorem is equivalent to the classical
 Levi-Civita's Theorem from
\cite{Levi-Civita}.

 \begin{Th}[Levi-Civita's Theorem]\hspace{-2mm}{\bf .} \label{LC}
 The following statements hold:
 \begin{enumerate}
 \item
 Let $L$ be a BM-structure on $(M^n, g)$.
 Let $x\in M^n$ be typical.  Then,  there exists
a coordinate system $\bar x=(\bar x_1,...,\bar x_m)$ in a
neighborhood  ${U}(x)$ containing $x$, where $\bar
x_i=(x_i^1,...,x_i^{k_i})$, $(1\le i\le m)$, such that $L$ is
diagonal
\begin{equation} \label{diagonal}
\diag(\underbrace{\phi_1,...,\phi_1}_{k_1},....,
\underbrace{\phi_m,...,\phi_m}_{k_m}),
\end{equation}
and the
 quadratic  form  of the
metric $g$  have the following form:
\begin{eqnarray}
g(\dot{\bar x}, \dot{\bar x})&=&\ P_1(\bar x)A_1(\bar x_1,\dot{\bar x}_1)+
\ P_2(\bar x)A_2(\bar x_2,\dot{\bar x}_2)+\cdots+\nonumber\\
&+&\ P_m(\bar x)A_m(\bar x_m,\dot{\bar x}_m), \label{g}
\end{eqnarray}

\noindent where $A_i(\bar x_i, \dot{\bar x}_i)$ are
 positive-definite quadratic forms  in  the velocities $\dot{\bar x}_i$
with coefficients depending  on $\bar x_i$,
\begin{eqnarray*}
P_i&\eqdef&(\phi_i-\phi_1)\cdots(\phi_{i}-\phi_{i-1})(\phi_{i+1}-\phi_i)\cdots
(\phi_m-\phi_i),
\end{eqnarray*}
and $\phi_1<\phi_2<...<\phi_m$   are smooth functions such that
$$
\phi_i=\left\{
\begin{array}{l}
\phi_i(\bar x_i),\mbox{\quad if}\quad k_i=1\\
{\rm constant},\quad\mbox{otherwise}.
\end{array}
\right.
$$
\item Let $g$ be a Riemannian metric and $L$ be a  (1,1)-tensor. If
in a neighborhood $U\subset M^n$ there exist coordinates $\bar
x=(\bar x_1,...,\bar x_m)$ such that $L$ and $g$ are given by
formulae (\ref{diagonal}, \ref{g}),
 then  the restriction of $L$ to $U$ is a BM-structure for the
 restriction of $g$ to $U$.
\end{enumerate}
\end{Th}

\begin{Rem}\hspace{-2mm}{\bf .} \label{eisen}
In Levi-Civita's coordinates from Theorem~\ref{LC}, the metric
$\bar g$ given by (\ref{bg}) has the form
\begin{equation}\label{rem2}\begin{array}{ccl}
\bar g(\dot{\bar x}, \dot{\bar x})&=&\rho_1P_1(\bar x)A_1(\bar
x_1,\dot{\bar x}_1)+ \rho_2P_2(\bar x)A_2(\bar x_2,\dot{\bar
x}_2)+\cdots+\nonumber\\ &+&\rho_mP_m(\bar x)A_m(\bar
x_m,\dot{\bar x}_m),\end{array}
\end{equation}
where
 \begin{eqnarray*}
\rho_i  & = & \frac{1}{\phi_1^{k_1}...\phi^{k_m}_m}\frac{1}{\phi_i}
\end{eqnarray*}
The metrics $g$ and $\bar g $ are affine equivalent (i.e. they have the same Levi-Civita connections) if and only if
all functions $\phi_i$ are constant.
\end{Rem}

Let $p$ be a typical point  with respect to the BM-structure $L$.
Fix $i\in 1,...,m$ and a small neighborhood $U$ of $p$. At every
point of $U$, consider the eigenspace $V_i$ with the eigenvalue
$\phi_i$. If the neighborhood is small enough, it contains only
typical points and  $V_i$ is a distribution. By Theorem~\ref{nijenhuis},
it is integrable. Denote by $M_i(p)$
the integral manifold containing $p$.

Levi-Civita's Theorem says that the eigenvalues $\phi_j$, $j\ne
i$, are constant on $M_i(p)$, and that the restriction of  $g$ to
$M_i(p)$ is proportional to the restriction of $g$ to $M_i(q)$, if
it is possible to connect $q$ and $p$ by a line orthogonal to
$M_i$ and containing only typical points. Actually, in view of \cite{Haantjes},  the first
observation follows already from Theorem~\ref{nijenhuis}. We will
need the second observation later and formulate it as

\begin{Cor}\hspace{-2mm}{\bf .}\label{corlc}
Let $L$ be a BM-structure for a  connected Riemannian manifold  $(M^n,g)$.  Suppose the
curve $\gamma:[0,1]\to M^n$ contains only typical points and is
orthogonal to $M_i(p) $ at every point $p\in
\textrm{Image}(\gamma)$. Let the multiplicity of the eigenvalue
$\phi_i$ at every point of the curve be  greater than one. Then,
the restriction of the metric to $M_i({\gamma(0)})$ is
proportional to the restriction of the metric to
$M_i({\gamma(1)})$. (i.e. there exists a diffeomorphism of  a
small neighborhood  $ U_i(\gamma(0)) \subset M_i({\gamma(0)})$ to
a small neighborhood  $U_i(\gamma(1))\subset M_i({\gamma(1)})$
taking the restriction of the metric $g$ to $M_i({\gamma(0)})$ to
a metrics proportional to the restriction of the metric $g$ to
$M_i({\gamma(1)})$).

\end{Cor}

\begin{Def}\hspace{-2mm}{\bf .}  \label{defadjusted} Let  $(M^n, g)$ be a Riemannian manifold.
We say that the metric $g$ has a {\bf warped decomposition}
 near $x\in M^n$, if a
neighborhood $U^n$ of $x$ can be split in the direct product of
disks $D^{k_0}\times ... \times D^{k_m}$, $k_0+...+k_m=n$, such that the metric $g$ has the form
 \begin{equation} \label{vk}
g_0+\sigma_1 g_1+\sigma_2 g_2+...+\sigma_m g_m,
\end{equation}
where the $i${\rm th}   metric  $g_i$ is a Riemannian metric on
the corresponding disk $D^{k_i}$,  and functions $\sigma_i$ are
functions on the disk $D^{k_0}$. The metric
\begin{equation} \label{adjusted}
g_0+\sigma_1 dy_1^2+\sigma_2 dy_2^2+...+\sigma_m dy^2_m
\end{equation}
on $D^{k_0}\times \mathbb{R}^m$ is called  the {\bf adjusted
metric}.
\end{Def}

We will always assume that $k_0$   is   at least $1$.
Adjusted metric has a very clear geometric sense.
Take a point
$p=(p_0,...,p_m)\in D^{k_0}\times ...\times  D^{k_m}.$
 At every disk $D^{k_i}$, \  $i=1,...,m$,
consider a geodesic segment $\gamma_{i}\in D^{k_i}$ passing
through $p_i$. Consider the product
 $$
M_A:=D^{k_0}\times \gamma_1\times\gamma_2\times ...
\times{\gamma_m}
 $$
 as a submanifold of $D^{k_0}\times ... \times  D^{k_m}.$
 As it easily  follows from Definition~\ref{defadjusted},
 \begin{itemize}
 \item $
M_A$ is a totally geodesic submanifold.
\item The restriction of the metric (\ref{vk}) to $M_A$ is
(isometric to) the adjusted metric.
\end{itemize}

Comparing formulae (\ref{g},\ref{vk}),  we see that if $L$ has at least one
simple eigenvalue at a typical point,
 Levi-Civita's
Theorem gives us a warped decomposition  near every typical point
of $M^n$: the metric $g_0$ collects all $P_iA_i$ from (\ref{g})
such that $\phi_i$ has multiplicity one,  the metrics
$g_1$,...,$g_m$ coincide with $A_j$ for multiple $\phi_j$, and
$\sigma_j=P_j$.

\begin{Def}[\cite{solodovnikov1,solodovnikov2.5,solodovnikov2.6}]\hspace{-2mm}{\bf .}
 Let $K$ be a constant. A metric $g$  is
called a {\bf $V(K)$-metric} near $x\in M^n$ ($n\ge 3$),
 if there exist coordinates
in a neighborhood of  $x$
 such that $g$ has the Levi-Civita form~(\ref{g}) such that
  the
 adjusted metric  has constant sectional curvature $K$.
\end{Def}

The  definition above is independent of
the choice of the presentation of $g$ in Levi-Civita's form:

\begin{Th}[\cite{solodovnikov1,solodovnikov2,solodovnikov2.5,solodovnikov2.6}] \hspace{-2mm}{\bf .} \label{maximal}
Suppose $g$  is a $V(K)$-metric near $x\in M^n$. Assume $n\ge 3$. The following statements hold:
\begin{enumerate}
\item If there exists another presentation   of $g$ (near $x$)   in the form
(\ref{g}), then  the sectional curvature of the adjusted metric
  constructed for this other decomposition is
constant and is equal to $K$.

\item  Consider the  metric (\ref{g}).   For every $i=1,...,m$, denote
\begin{equation}\label{Ki}
\frac{g(grad(P_i), grad(P_i))}{4\, P_i} + K P_i\end{equation}   by
$K_i$.
Then,
 the metric (\ref{vk}) has
       constant sectional curvature if and only  if for every  $i\in 0,...,m$ such that
$k_i>1$ the metric  $A_i$
       has constant sectional  curvature $K_i$. More precisely, if the metric
       (\ref{vk}) is a $V(K)$-metric,  if $k_1>1$ and if the metric $A_1$  has
       constant sectional  curvature $K_1$, then the metric $g_0+P_1 A_1$  has constant sectional curvature $K$.

\item For a fixed presentation of $g$ in the Levi-Civita form (\ref{g}),
for every $i$ such that $k_i>1$,    $ K_i$ is a constant.
\end{enumerate}
\end{Th}

The first statement of Theorem~\ref{maximal} is proved in \S 3, \S 7 of
\cite{solodovnikov1}.  In the form sufficient for our paper, it
  appeared already in
\cite{Vries}; although it is  hidden there.
The second and the third   statements    can
 be found,
 for example,  in \S 8 of
 \cite{solodovnikov1}.
 The relation between
$V(K)$-metrics and BM-structures is given by

\begin{Th}[\cite{solodovnikov2.5,solodovnikov2.6}]\hspace{-2mm}{\bf .} \label{solodovnikov}
Let $(D^n,g)$ be a disk of dimension $n\ge 3$  with two
BM-structures $L_1$ and $L_2$ such that every point of the  disk
is typical with respect to both structures and the BM-structures
$\Id, \ L_1, \ L_2$ are linearly independent.  Assume that at least one eigenvalue of $L_1$ is simple.

 Then, $g$ is a
$V(K)$-metric near every point.
\end{Th}
A partial case of this theorem is
\begin{Th}[Fubini's Theorem \cite{Fubini1,Fubini2}]\hspace{-2mm}{\bf .} \label{fubini}\label{constant}
Let $(D^n,g)$ be a disk of dimension $n\ge 3$  with two
BM-structures $L_1$ and $L_2$ such that $N_{L_1}=N_{L_2}=n$ at
every point.  If the BM-structures \ $\Id, \ L_1, \ L_2$ are
linearly independent,  then $g$ has constant curvature.
\end{Th}

\section{Global theory of projectively equivalent metrics}
\label{global}

The goal of this section is to provide necessary tools for
Section~\ref{proof}.

\subsection{The eigenvalues of a BM-structure  are globally
ordered} \label{proofordered}

\noindent Within this section we assume that $(M^n,g)$ is a
Riemannian manifold such that every two points can be connected by
a geodesic.  Our goal is to prove Corollary~\ref{ordered12}. We
need the following technical lemma. For every fixed
$v=(\xi_1,\xi_2, ... ,\xi_n)\in T_xM^n$,
 the function (\ref{integral}) is
 a  polynomial  in $t$.
Consider the roots of this  polynomial.
 From the proof of Lemma~\ref{technical},
it will be clear that they are real. We denote them by
 $$
t_1(x,v)\le t_2(x,v)\le ... \le  t_{n-1}(x,v). $$
 \begin{Lemma}[\cite{hyperbolic,dubrovin}]\hspace{-2mm}{\bf .} \label{technical} The following holds
 for every  $i\in \{1,...,n-1\}$:
\begin{enumerate}
 \item For every  $v \in T_x M^n$,
  $$
    \lambda_i(x)\le t_i(x,v) \le \lambda_{i+1}(x).
  $$
     In particular,  if $\lambda_i(x) = \lambda_{i+1}(x)$, then
$t_i(x,v)=\lambda_i(x) = \lambda_{i+1}(x)$.

  \item If $\lambda_i(x) < \lambda_{i+1}(x)$, then for every   $\tau\in \mathbb R $
the
   Lebesgue
measure of the set
 $$
V_\tau\subset   T_x M^n, \ \ V_\tau\eqdef \{ v\in T_xM^n: t_i(x,v)
= \tau\}, $$
  is zero.
\end{enumerate}
 \end{Lemma}

\noindent{\bf  Proof:} By definition, the
tensor $L$ is self-adjoint with respect to $g$. Then, for every
$x\in M^n$,
 there exist  "diagonal" coordinates  in $T_xM^n$ such that   the
 metric $g$ is given by the    diagonal matrix
 $\diag(1,1, ... ,1)$  and
the tensor  $L$ is given by the
  diagonal  matrix $\diag(\lambda_1,\lambda_2, ... ,\lambda_n)$.
Then, the tensor
 (\ref{st}) reads:
\begin{eqnarray*}
S_t& = & \det(L-t\Id)(L-t\Id)^{(-1)} \\ &=& \diag(\Pi_1(t),
\Pi_2(t), ... ,\Pi_n(t)),
\end{eqnarray*}
where the polynomials  $\Pi_i(t)$ are  given by the formula $$
\Pi_i(t)\eqdef (\lambda_1-t)(\lambda_2-t) ...
(\lambda_{i-1}-t)(\lambda_{i+1}-t) ...
(\lambda_{n-1}-t)(\lambda_n-t). $$ Hence, for every
$v=(\xi_1,...,\xi_n)\in T_xM^n$, the polynomial $I_t(x,v)$ is
given by
\begin{equation}\label{33}
I_t=\xi_1^2\Pi_1(t)+ \xi_2^2\Pi_2(t)+ ... +\xi_n^2\Pi_n(t).
\end{equation}
Evidently,  the coefficients of the polynomial $I_t$ depend
continuously on the eigenvalues $\lambda_i$ and on the components
 $\xi_i$.
Then, it is sufficient to prove the first statement of the lemma
assuming that the eigenvalues  $\lambda_i$ are all different and
that  $\xi_i$
 are
non-zero.  For every  $\alpha\ne i$,
 we evidently
have $\Pi_\alpha(\lambda_i)\equiv 0$. Then, $$ I_{\lambda_i}
=\sum_{\alpha=1}^{n}\Pi_\alpha (\lambda_i) \xi_\alpha^2 =
\Pi_i(\lambda_i) \xi_i^2. $$ Hence  $I_{\lambda_i}(x,v)$ and
$I_{\lambda_{i+1}}(x,v)$ have different signs.
 Hence,
 the open interval $]\lambda_i,\lambda_{i+1}[$ contains
a root of the polynomial $I_t(x,v)$. The degree of the polynomial
$I_t$ is equal $n-1$; we have
 $n-1$ disjoint intervals; every
 interval contains at least one
root so that all roots are real and the  $i$th  root
 lies between
$\lambda_i$ and $\lambda_{i+1}$.
  The first statement of  Lemma~\ref{technical} is proved.

Let us prove the second statement of Lemma~\ref{technical}.
Assume  $\lambda_i<\lambda_{i+1}$. Suppose  first $\lambda_i<\tau<
\lambda_{i+1}$. Then, the set $$  V_\tau\eqdef \{ v\in T_xM^n:
t_i(x,v) = \tau\},$$ consists of the points $v$ such that  the function
$I_{\tau}(x,v)\eqdef (I_t(x,v))_{|t= \tau}$ is zero.  Then,    $V_{\tau}$ is
  a nontrivial quadric in $T_xM^n\equiv \mathbb{R}^n$ and, hence, has
 zero measure.

Now suppose  $\tau$ is  an  endpoint of the interval $[\lambda_i,
\lambda_{i+1}]$. Without loss of generality, we can assume $\tau
= \lambda_i$. Let $k$ be the multiplicity of the eigenvalue
$\lambda_i$. Then, every coefficient $\Pi_\alpha(t)$ of the
quadratic form
 (\ref{33}) has  the  factor
$(\lambda_i - t)^{k-1}$. Hence, $$ \hat I_t \eqdef
\frac{I_t}{(\lambda_i - t)^{k-1}} $$ is a polynomial in $t$ and
$\hat I_\tau$ is a nontrivial quadratic form. Evidently, for every
point   $v \in V_\tau$, we have $\hat I_\tau(v) = 0$ so that the
set $V_\tau$ is a subset of a nontrivial quadric in
 $T_xM^n$ and, hence, has   zero  measure.
  Lemma~\ref{technical}   is
proved.

{\noindent \bf Proof of Corollary~\ref{ordered12}:} The first
statement of Corollary~\ref{ordered12} follows immediately from
the first statement of Lemma~\ref{technical}: Let us join the
points $x,y\in M^n$  by a geodesic $\gamma:\mathbb{R}\to M^n$,
$\gamma(0)=x, \ \gamma(1)=y$. Consider the one-parametric family
of integrals
$
I_t(x,v)
$
and the  roots $$ t_1(x,v) \le t_2(x,v) \le ... \le t_{n-1}(x,v).
$$

 By Theorem~\ref{integrability},
 every   root $t_i$ is   constant
on every  orbit $(\gamma,\dot\gamma)$ of the geodesic flow of $g$
so that $$ t_i(\gamma(0),\dot\gamma(0)) =
t_i(\gamma(1),\dot\gamma(1)). $$ Using Lemma~\ref{technical}, we
obtain $$ \lambda_i(\gamma(0))\le t_i(\gamma(0),\dot\gamma(0)), \
\ \ \mbox{and}\ \ \ \
 t_i(\gamma(1),\dot\gamma(1)) \le \lambda_{i+1}(\gamma(1)).
$$ Thus $\lambda_i(\gamma(0))\le \lambda_{i+1}(\gamma(1))$ and the
first statement of Corollary~\ref{ordered12} is proved.

Let us prove the second statement of Corollary~\ref{ordered12}.
 Suppose
 $\lambda_i(y) = \lambda_{i+1}(y)$
 for every point $y$ of  some  subset
$V\subset M^n$. Then,  $\lambda_i$ is a constant on $V$
(i.e.  $\lambda_i$ is independent
 of $y\in V$). Indeed, by the first statement
 of Corollary~\ref{ordered12},
$$ \lambda_i(y_0)\le \lambda_{i+1}(y_1)   \  \ \mbox{and} \  \
\lambda_i(y_1) \le  \lambda_{i+1}(y_0), $$ so that $
\lambda_i(y_0)=\lambda_i(y_1)=\lambda_{i+1}(y_1)=\lambda_{i+1}(y_0)$
for every  $y_0, y_1\in V$.

We denote this constant by $\tau$.
 Let us join the point $x$
with every point of $V$ by all possible geodesics. Consider the
set $V_\tau \subset T_xM^n$ of the initial velocity vectors (at
the point $x$)
 of these geodesics.

 By the first statement of   Lemma~\ref{technical},
 for every  geodesic $\gamma$ passing
 through at least one  point
of $V$,  the value $t_i(\gamma,\dot\gamma)$
 is equal to $\tau$.
By the second  statement of   Lemma~\ref{technical}, the measure
of  $V_\tau $ is zero. Since the set $V$ lies in the image
of the exponential mapping of  $V_\tau$,
 the measure of the set $V$ is also zero. Corollary~\ref{ordered12}  is
proved.

\subsection{Local theory: behavior of BM-structure near non-typical points. }
Within this section  we assume that $L$ is a BM-structure on  a
connected $(M^n,g)$. As in Section~\ref{integrals}, we denote by
$\lambda_1(x)\le ... \le \lambda_n(x)$ the eigenvalues of $L$, and
by $N_L(x)$ the number of different eigenvalues of $L$ at $x\in
M^n$.

\begin{Th}\hspace{-2mm}{\bf .} \label{isol}
 Suppose the eigenvalue $\lambda_1$ is not constant,
the eigenvalue $\lambda_2$ is constant
  and
 $N_L=2$ in a typical
point. Let $p$ be a non-typical point.  Then, the following
statements hold:
\begin{enumerate}
\item The spheres of small radius with center in  $p$ are orthogonal to the
eigenvector of $L$ corresponding to $\lambda_1$, and tangent to
the eigenspace of $L$ corresponding to $\lambda_2$. In particular,
the points $q$ such that  $\lambda_1(q)=\lambda_2$ are isolated.
\item For every sphere of small radius with center in  $p$,
the restriction of $g$ to the sphere   has constant sectional
curvature.
\end{enumerate}
\end{Th}
{\bf Proof: }  Since $\lambda_1$ is not constant, it is  a simple
eigenvalue in every typical point. Since $N_L=2$, the roots
$\lambda_2,\lambda_3,...,\lambda_n$ coincide at every point and
are constant. We denote this constant by $\lambda$.  By
Lemma~\ref{technical}, at every point $(x,\xi)\in T_xM^n$, the
number  $\lambda$ is a root of multiplicity at least $n-2$ of the
polynomial $I_t(x,\xi)$. Then, $$
I'_t(x,\xi):=\frac{I_t(x,\xi)}{(\lambda-t)^{n-2}}$$ is a linear
function in $t$ and, for every fixed $t$, is an integral of the
geodesic flow of $g$. Denote by $\tilde I:TM\to \mathbb{R}$ the
function $$ \tilde I(x,\xi):= I'_{\lambda}(x,\xi):=\left(
I'_t(x,\xi)\right)_{|t=\lambda}.$$ Since $\lambda$ is a constant,
the function $\tilde I$ is an integral of the geodesic flow of
$g$. At every tangent space $T_xM^n$,  consider the coordinates
such that the metric is given by $\diag(1,...,1)$ and $L$ is given
by $\diag(\lambda_1,\lambda,...,\lambda)$. By direct calculations
we see that the restriction of $\tilde I$ to $T_xM^n$ is given by
(we assume $\xi=(\xi_1,\xi_2,...,\xi_n)$)
 $$ \tilde I_{|T_xM^n}(\xi) =
(\lambda_1(x)-\lambda)(\xi_2^2+...+\xi_n^2). $$ Thus, for every
geodesic $\gamma$ passing through $p$,  the value of $\tilde
I(\gamma(\tau),\dot \gamma(\tau))$ is zero.  Then, for every
typical point of such geodesic, since $\lambda_1<\lambda$,  the
components $\xi_2,...,\xi_n$ of  the  velocity vector vanish.
Finally,  the velocity vector is an eigenvector of $L$ with the
eigenvalue $\lambda_1$.

Then, the points where  $\lambda_1=\lambda$ are isolated:
otherwise we can pick two such points $p_1$ and $p_2$ lying in a
ball with radius less than the radius of injectivity. Then, for
almost every point $q$  of the ball, the geodesics connecting $q$ with
$p_1$ and $p_2$ intersect transversally at
$q$. Then, the point $q$ is non-typical; otherwise the eigenspace
of $\lambda_1$ contains the velocity vectors of geodesics and is
not  one-dimensional. Finally, almost every point of the ball
is not typical, which contradicts Corollary~\ref{typical}. Thus,
the points where  $\lambda_1=\lambda$ are isolated.

It is known (Lemma of Gau\ss{}),   that the geodesics passing
through $p$ intersect
 the
spheres of small radius with center in  $p$ orthogonally. Since
the velocity vectors of such geodesics are eigenvectors of $L$
with eigenvalue $\lambda_1$,  then the eigenvector with eigenvalue
$\lambda_1$ is orthogonal to the spheres of small radius with
center in $p$. Since $L$  is self-adjoint, the spheres are tangent
to the eigenspaces  of $\lambda$. The first statement of
Theorem~\ref{isol} is proved.

The second statement of Theorem~\ref{isol} is trivial, if $n=2$.
In order to prove the second statement for $n\ge 3$, we will use
Corollary~\ref{corlc}. The role of the curve $\gamma$ from
Corollary~\ref{corlc} plays the geodesic passing through $p$. We
put $i=2$. By the first statement  of Theorem~\ref{isol}, \
$M_i(x)$ are spheres with center in $p$. Then, by
Corollary~\ref{corlc}, for every sufficiently small spheres
$S_{\epsilon_1}$ and $S_{\epsilon_2}$ with center in $p$,  \ the
restriction of $g$ to the first sphere is proportional to the
restriction of $g$ to the second sphere. Since for very small
values of $\epsilon$ the metric in a $\epsilon$-ball is very close to the
Euclidean metric, the  restriction of $g$ to the $\epsilon$-sphere
is close to the round metric of the sphere. Thus, the restriction
of $g$ to every (sufficiently small) sphere with center in $p$  has constant sectional
curvature. Theorem~\ref{isol} is proved.

\begin{Th}\hspace{-2mm}{\bf .} \label{add}
Suppose $N_L=3$ at a typical point and there exists a point where
$N_L=1$.  Then,  the following statements hold:

\begin{enumerate} \item There exist points $p_1$, $p_n$  such that
$\lambda_1(p_1)<\lambda_2(p_1)=\lambda_n(p_1)$ and
$\lambda_1(p_n)=\lambda_2(p_n)<\lambda_n(p_n)$.

\item The points $p$ such that  $N_L(p)=1$ are isolated.
\end{enumerate}
\end{Th}
{\bf Proof:}  Let us prove  the first statement.  Suppose
$\lambda_1(p_2)=\lambda_2(p_2)=...=\lambda_n(p_2)$ and the number
of different eigenvalues of $L$ at a typical point equals three.
Then, by Corollary~\ref{ordered12}, the eigenvalues
$\lambda_2=...=\lambda_{n-1}$ are constant. We denote this
constant by $\lambda$. Take a ball $B$ of small radius with center
in $p_2$. We will prove that this ball has a point $p_1$ such that
$\lambda_1(p_1)<\lambda_2=\lambda_n(p_1)$; the proof that there
exists a point where $\lambda_1=\lambda_2<\lambda_n$ is similar.
Take $p\in B$ such that $\lambda_1(p)<\lambda$ and $\lambda_1(p)$
is a regular value of the function $\lambda_1$. Denote by $\check
M_1(p)$ the connected component of $\{q\in M^n:\
\lambda_1(q)=\lambda_1(p)\}$ containing the point $p$. Since
$\lambda_1(p)$ is  a regular value, $\check M_1(p)$ is a
submanifold of codimension 1. Then, there exists a point $p_1\in
\check M_1(p)$ such that the distance from this point to $p_2$ is
minimal over all points of $\check M_1(p)$.

Let us  show that $\lambda_1(p_1)<\lambda=\lambda_n(p_1)$. The
inequality $ \lambda_1(p_1)<\lambda$ is fulfilled by definition,
since $p_1\in \check M_1(p)$. Let us  prove that
$\lambda_n(p_n)=\lambda$.

 Consider the shortest geodesic $\gamma$
connecting $p_2$ and $p_1$. We will assume $\gamma(0)=p_1$ and
$\gamma(1)=p_2$. Consider the values of the roots $t_1\le ... \le
t_{n-1}$ of the polynomial $I_t$ at the  points of  the geodesic orbit
$(\gamma,\dot \gamma)$. Since $I_t$ are integrals, the roots $t_i$
are independent  of the point of the orbit. Since the geodesic
passes  through the point where $\lambda_1=...=\lambda_n$, by
Lemma~\ref{technical}, we have
\begin{equation} \label{ert}
 t_1= ...
=t_{n-1}=\lambda. \end{equation}
 Since the distance from
$p_1$ to $p_2$ is minimal over all points of $\check M_1$, the
velocity vector $\dot\gamma(0)$ is orthogonal to $\check M_1$. By
 Theorem~\ref{nijenhuis} and \cite{Haantjes}, the sum of eigenspaces of $L$
corresponding to $\lambda$ and $\lambda_n$ is  tangent to $\check
M_1$. Hence,  the vector $\dot\gamma(0)$ is  an eigenvector of $L$
with eigenvalue $\lambda_1$.

 At the tangent space
$T_{p_1}M^n$, choose a coordinate system such that $L$  is
diagonal $\diag(\lambda_1,...,\lambda_n)$ and $g$ is Euclidean
$\diag(1,...,1)$. In this coordinate system, $I_t(\xi)$ is given
by (we assume $\xi=(\xi_1,...,\xi_n)$) \begin{equation}\label{tmp1}
(\lambda-t)^{n-3}\left((\lambda_n-t)(\lambda-t)\xi_1^2+(\lambda_{n}-t)(\lambda_1-t)(\xi_2^2+...+\xi_{n-1}^2)+
(\lambda_1-t)(\lambda-t)\xi_{n}^2\right). \end{equation}  Since $\dot\gamma(0)$
is an eigenvector of $L$ with eigenvalue $\lambda_1$, the last
$n-1$ components of $\dot\gamma(0)$ vanish, so that
$t_{n-1}=\lambda_n$. Comparing this with (\ref{ert}), we see that
$\lambda_n(p_1)=\lambda$. The first statement of
Theorem~\ref{add}  is proved.

Now let us prove the second statement. We suppose that in a small convex  ball
$B\subset M^n$ there exists four  points $p', p'',p''',p''''$
with $N_L=1$, and will find
 a contradiction.    By Corollary~\ref{typical}, almost every point
 $p$ of the
ball is typical. Clearly, for almost every   typical  point $p$
 of the ball  the geodesics connecting the point  with
 $p',  p'', p''',p''''$ intersect
mutually-transversally in $p$. Since these geodesic pass through
 points  where $\lambda_1=...=\lambda_n$, by Lemma~\ref{technical},  the roots
$t_1,t_2,...,t_{n-1}$ on the geodesics are all  equal to $\lambda$.
If  the point $p$ is typical, the restriction of  $I_t$ to
$T_pM^n$  has the form (\ref{tmp1}). Then, if
$(\xi_1',...,\xi_n')$,  $(\xi_1'',...,\xi_n'')$,  $(\xi_1''',...,\xi_n''')$,
$(\xi_1'''',...,\xi_n'''')$
 are the coordinates of velocity vectors of the geodesics at $p$,
the sums $((\xi_2')^2+...+(\xi_{n-1}')^2)$,
$((\xi_2'')^2+...+(\xi_{n-1}'')^2)$,
$((\xi_2''')^2+...+(\xi_{n-1}''')^2)$,
$((\xi_2'''')^2+...+(\xi_{n-1}'''')^2)$
vanish and  the following system of equations holds:
$$
\left\{\begin{array}{ccc}
(\lambda_n-\lambda)(\xi_1')^2+
(\lambda_1-\lambda)(\xi_{n}')^2&=&0\cr
(\lambda_n-\lambda)(\xi_1'')^2+
(\lambda_1-\lambda)(\xi_{n}'')^2&=&0\cr
(\lambda_n-\lambda)(\xi_1''')^2+
(\lambda_1-\lambda)(\xi_{n}''')^2&=&0\cr
(\lambda_n-\lambda)(\xi_1'''')^2+
(\lambda_1-\lambda)(\xi_{n}'''')^2&=&0
\end{array}
\right.
$$
Thus, at  least two of  the vectors
$(\xi_1',...,\xi_n')$ and $(\xi_1'',...,\xi_n'')$ are proportional. Then, there exists a pair of geodesics that   do not intersect transversally at $p$.
The contradiction shows that
 the points where
$N_L=1$ are isolated. Theorem~\ref{add} is proved.

\subsection{Splitting Lemma} \label{splitting}

\begin{Def}\hspace{-2mm}{\bf .} \label{lpr}
A {\bf local-product structure} on $M^n$ is the triple $(h, B_r,
B_{n-r})$, where  $h$ is a Riemannian metric and $B_r$, $B_{n-r}$
are transversal foliations of dimensions   $r$ and $n-r$,
respectively (it is assumed that $1\le r< n$),
 such that every    point $p\in M^n$
has  a   neighborhood $U´(p)$   with coordinates $$ (\bar x,\bar
y)= \bigr((x_1,x_2,...x_r),(y_{r+1},y_{r+2},...,y_n)\bigl) $$ such
that  the  $x$-coordinates are constant on every leaf  of the
foliation $B_{n-r}\cap U´(p)$, the $y$-coordinates are constant on
every leaf of the foliation $B_{r}\cap U´(p)$,   and
 the metric  $h$ is block-diagonal   such that the first
($r\times r$) block depends on the $x$-coordinates and the last
$((n-r)\times (n-r))$ block depends on the $y$-coordinates.
\end{Def}

A model example of manifolds with local-product structure is the
direct product of two Riemannian manifolds $(M_1^r, g_1)$ and
$(M_2^{n-r}, g_2)$.  In this case, the leaves  of the foliation
$B_r$ are the products  of $M_1^r$ and  the points of $M_2^{n-r}$,
the leaves  of the foliation $B_{n-r}$ are the products  of  the
points of $M_1^{r}$ and  $M_2^{n-r}$, and the  metric $h$ is the
product metric $g_1+g_2$.

Below we assume that

(a) $L$ is a BM-structure for a connected   $(M^n,g)$.

(b) There exists $r$, $1\le r<n$,  such  that
$\lambda_r<\lambda_{r+1}$ at every  point of $M^n$.


We will show that (under the assumptions (a,b)) we can naturally
define a local-product structure $(h, B_r, B_{n-r})$
 such that the (tangent spaces to) leaves of $B_r$ and $B_{n-r}$ are
 invariant with respect to $L$, and such that the restrictions
$L_{|B_r}$, $L_{|B_{n-r}}$ are BM-structures for the metrics
 $h_{|B_r}$, $h_{|B_{n-r}}$, respectively.

 At every
 point $x\in M^n$, denote by $V^r_x$
  the  subspaces  of $T_xM^n$ spanned by the eigenvectors of
$L$ corresponding to the eigenvalues $\lambda_1,...,\lambda_r$.
Similarly,  denote by $V^{n-r}_x$
  the  subspaces  of $T_xM^n$ spanned by the eigenvectors of
$L$ corresponding to the eigenvalues
$\lambda_{r+1},...,\lambda_n$. By assumption, for every  $i,j$
such that
 $i\le r < j$, we have $\lambda_i\ne \lambda_j$ so that
$V^r_x$ and $V^{n-r}_x$ are  two
 smooth distributions  on $M^n$. By Theorem~\ref{nijenhuis},
the distributions are integrable so that they define two
transversal     foliations $B_r$ and $B_{n-r}$ of dimensions $r$
and $n-r$, respectively.

By construction,   the distributions $V_r$ and
 $V_{n-r}$ are invariant with respect to $L$.  Let us denote by
$L_r$,  $L_{n-r}$ the restrictions of $L$ to $V_r$ and
 $V_{n-r}$, respectively.
We will denote by $\chi_r$, $\chi_{n-r}$  the characteristic
polynomials of $L_r$,  $L_{n-r}$, respectively.
 Consider the
 (1,1)-tensor
$$ C\eqdef  \left((-1)^{r}\chi_r(L)+\chi_{n-r}(L)\right) $$ and
the metric $h$ given by the relation $$ h(u,v)\eqdef
g(C^{-1}(u),v) $$ for every   vectors $u,v$. (In the
 tensor notations, the metrics  $h$ and $g$ are related by
$g_{ij}=h_{i\alpha}C_j^\alpha$).

\begin{Lemma}[Splitting Lemma]\hspace{-2mm}{\bf .}
 \label{local-direct-product} The following statements hold:

 \begin{enumerate}
  \item The triple $(h, B_r, B_{n-r})$
 is  a local-product structure  on $M^n$.

  \item   For every  leaf of $B_r$, the restriction  of $L$ to it is a BM-structure for the restriction of
   $h$  to it. For every  leaf of $B_{n-r}$, the restriction of $L$ to it is a BM-structure for the restriction of
   $h$  to it.
\end{enumerate}
\end{Lemma}

\noindent {\bf Proof:}  First of all,  $h$ is a well-defined
 Riemannian metric.
Indeed, take an arbitrary point $x\in M^n$. At the tangent space
to this point,
  we can find a coordinate system such that  the tensor $L$  and the metric $g$
are diagonal.  In this coordinate system, the characteristic
polynomials  $\chi_r$, $\chi_{n-r}$ are given by
\begin{equation}\label{lcr}
\begin{array}{lcr}
(-1)^{r}\chi_r&=& (t-\lambda_1)(t-\lambda_2)...(t-\lambda_r)   \\
\chi_{n-r}&=& (\lambda_{r+1}-t)(\lambda_{r+2}-t)...(\lambda_n-t).
\end{array}
\end{equation}
Then,  the (1,1)-tensor \ \
$
C=\left((-1)^{r}\chi_r(L)+\chi_{n-r}(L)\right)$ \ \
\noindent  is   given by the  diagonal matrices
\begin{eqnarray}
\hspace{-5mm}\diag\left(\prod_{j=r+1}^n(\lambda_{j}-\lambda_1),...,
      \prod_{j=r+1}^n(\lambda_{j}-\lambda_r),
 \prod_{j=1}^r(\lambda_{r+1}-\lambda_j),...,
      \prod_{j=1}^r(\lambda_{n}-\lambda_j)\right), \label{chi1}
\end{eqnarray}
We see that the tensor is    diagonal and that
 all   diagonal
components
  are positive.
Then,  the tensor  $C^{-1}$ is   well-defined and $h$ is  a
Riemannian metric.

By construction, $B_r$ and $B_{n-r}$ are well-defined transversal
foliations of  dimensions $r$ and $n-r$.  In order to prove
Lemma~\ref{local-direct-product}, we need to verify
 that,  locally,  the triple  $(h,B_r,B_{n-r})$
  is  as in Definition~\ref{lpr},  that the
restriction of $L$   to a leaf   is  a BM-structure for the
restriction of $h$ to the  leaf.

It is  sufficient to verify these  two  statements
 at almost every point of $M^n$. Indeed,
  since  the    foliations and the  metric are
globally-given and smooth, if the  restrictions of $g$ and  $L$ to a leaf
 satisfies Definitions~\ref{lpr},\ref{bm} at almost every point, then it
 satisfies Definitions~\ref{lpr},\ref{bm}  at every point.

Thus, by
 Corollary~\ref{typical}, it is sufficient to prove the Splitting Lemma near
 every typical point.  Consider
 Levi-Civita's coordinates
  $\bar x_1,...,\bar  x_m$ from Theorem~\ref{LC} near a typical point. As in
  Levi-Civita's Theorem, we denote by $\phi_1<...<\phi_m$ the
  different eigenvalues of  $L$. In Levi-Civita's  coordinates,  the matrix of $L$ is diagonal
  $$
\diag\left(\underbrace{\phi_1,...,\phi_1}_{k_1},...,\underbrace{\phi_m,...,\phi_m}_{k_m}\right)
=\diag(\lambda_1,...,\lambda_n).
  $$

Consider $s$ such that $\phi_s=\lambda_r$ (clearly,
$k_1+...+k_s=r$). Then, by constructions of the foliations $B_r$
and $B_{n-r}$,
 the coordinates
$\bar x_1,...,\bar x_s$
  are  constant on every leaf  of the foliation $B_{n-r}$,
the coordinates $\bar x_{s+1},...,\bar x_m$
 are
constant on every leaf of the foliation $B_{r}$. The coordinates
$\bar x_1,...,\bar x_s$ will play the role of $x$-coordinates from
Definition~\ref{lpr}, and the coordinates $\bar x_{s+1},...,\bar
x_m$ will play the role of $y$-coordinates from
Definition~\ref{lpr}.

Using (\ref{chi1}), we see that, in Levi-Civita's coordinates, $C$
is given by
\begin{eqnarray*}
\diag\left(\underbrace{\prod_{j=s+1}^m(\phi_{j}-\phi_1)^{k_j},...,\prod_{j=s+1}^m(\phi_{j}-\phi_1)^{k_j}}_{k_1},...,
      \underbrace{\prod_{j=s+1}^m(\phi_{j}-\phi_s)^{k_j},...,\prod_{j=s+1}^m(\phi_{j}-\phi_s)^{k_j}}_{k_s},\right.\\\left.
 \underbrace{\prod_{j=1}^s(\phi_{s+1}-\phi_j)^{k_j},...,\prod_{j=1}^s(\phi_{s+1}-\phi_j)^{k_j}}_{k_{s+1}},...,
    \underbrace{\prod_{j=1}^s(\phi_{m}-\phi_j)^{k_j},...,\prod_{j=1}^s(\phi_{m}-\phi_j)^{k_j}}_{k_{m}}
      \right).
\end{eqnarray*}

Thus,
 $h$ is  given by

\begin{equation}
\begin{array}{ccccc}
 h(\dot{\bar x}, \dot{\bar x})&=& \tilde
P_1A_1(\bar x_1,\dot{\bar x}_1)&+...+& \tilde  P_s A_s( \bar x_s,
\dot{\bar x}_s)  \label{newg1}
\\ &+&\tilde P_{s+1}A_{s+1}(\bar x_{s+1},\dot{\bar
x}_{s+1})&+...+&\tilde P_{m}A_{m}(\bar x_m,\dot{\bar x}_m),
\label{newg2}
\end{array}
\end{equation}

where the functions $\tilde P_i$ are as follows: for $i\le s$,
they  are  given by
\begin{eqnarray*}
\tilde P_i &\eqdef &
{(\phi_i-\phi_1)...(\phi_i-\phi_{i-1})(\phi_{i+1}-\phi_i)...(\phi_s-\phi_i)}
\prod_{\begin{array}{c}j=1\\ j\ne
i\end{array}}^s|\phi_i-\phi_j|^{1-k_j}.
\end{eqnarray*}
For   $i> s$,   the functions $\tilde P_i$ are given by
\begin{eqnarray*}
\tilde P_i & \eqdef &
{(\phi_i-\phi_{s+1})...(\phi_i-\phi_{i-1})(\phi_{i+1}-\phi_i)...(\phi_m-\phi_i)}
\prod_{\begin{array}{c}j=s+1\\ j\ne
i\end{array}}^m|\phi_i-\phi_j|^{1-k_j} .
 \end{eqnarray*}

Clearly, $|\phi_i-\phi_j|^{k_j-1}$ can  depend on the variables
$\bar x_i$ only. Then, the products
$$\prod_{\begin{array}{c}j=1\\ j\ne
i\end{array}}^s|\phi_i-\phi_j|^{1-k_j} \ ,
 \ \ \ \prod_{\begin{array}{c}j=s+1\\ j\ne
i\end{array}}^m|\phi_i-\phi_j|^{1-k_j}$$ can be hidden in $A_i$,
i.e. instead of $A_i$ we consider
$$\tilde A_i  \eqdef  \prod_{\begin{array}{c}j=1\\ j\ne
i\end{array}}^s|\phi_i-\phi_j|^{1-k_j} \ A_i \ \ \  \textrm{ for $i\le s$ and}
$$
$$\tilde A_i  \eqdef \prod_{\begin{array}{c}j=s+1\\ j\ne
i\end{array}}^m|\phi_i-\phi_j|^{1-k_j} \ A_i \ \  \ \textrm{for $i>s$.}$$

 Finally,
  the restriction of the metric to the leaves of $B_r$
has  the form from Levi-Civita's Theorem.  Hence, the
restriction of $L$ is a BM-structure for it. We see that the
leaves of $B_r$ are orthogonal to leaves of $B_{n-r}$, and that
the restriction of $h$ to $B_r$  ($B_{n-r}$, respectively) is
precisely the first row of (\ref{newg1}) (second row of
(\ref{newg2}), respectively) and depends on the coordinates $\bar
x_1,..., \bar x_s$ ($\bar x_{s+1},..., \bar x_m$, respectively)
only.
  Lemma~\ref{local-direct-product}
is proved.

 Let $p$ be a typical point  with respect to
the BM-structure $L$. Fix $i\in 1,...,n$. At every point of $M^n$,
consider the eigenspace $V_i$ with the eigenvalue $\lambda_i$.
$V_i$ is a distribution near $p$. Denote by $M_i(p)$ its integral
manifold containing $p$.

\begin{Rem}\hspace{-2mm}{\bf .}  \label{rest} The following statements  hold:
\begin{enumerate}
\item If $\lambda_i(p)$ is multiple, the restriction of $g$ to $M_i(p)$ is proportional to the
restriction of $h$ to $M_i(p)$.
\item The restriction of $L$ to $B_r$ does not depend on the
coordinates $y_{r+1},...,y_n$ (which are  coordinates $\bar
x_{s+1}, ... , \bar x_m$ in the notations in proof of
Lemma~\ref{local-direct-product}). The restriction of $L$ to
$B_{n-r}$ does not depend on the coordinates $x_{1},...,x_r$
(which are coordinates $\bar x_{1}, ... , \bar x_s$ in the
notations in proof of Lemma~\ref{local-direct-product}).
\end{enumerate}
\end{Rem}

   Combining
Lemma~\ref{local-direct-product} with Theorem~\ref{isol}, we
obtain
\begin{Cor}\hspace{-2mm}{\bf .}\label{corisol}
Let $L$ be BM-structure on  connected  $(M^n,g)$. Suppose there
 exist $i\in
1,...,n$ and $p\in M^n$ such that:
\begin{itemize}
\item $\lambda_i$ is multiple (with multiplicity $k\ge 2$)  at a typical point.
\item $\lambda_{i-1}(p)=\lambda_{i+k-1}(p)<\lambda_{i+k}(p)$,
\item The eigenvalue $\lambda_{i-1}$ is not constant.
\end{itemize}
Then, for every  typical point $q\in M^n$ which is  sufficiently
close to $p$,  $M_i(q)$ is diffeomorphic to the  sphere and the
restriction of $g$ to $M_i(q)$ has constant sectional curvature.
\end{Cor}

Indeed, take a small neighborhood of $p$ and  apply Splitting
Lemma~\ref{local-direct-product} two times: for $r=i+k-1$ and for
$r=i-2$. We obtain a metric $h$ such that locally, near $p$, the
manifold with this metric  is the Riemannian product of three
disks with  BM-structures, and BM-structure is the direct  sum of
these BM-structures. The second component of such decomposition
satisfies  the assumption of Theorem~\ref{isol}; applying
Theorem~\ref{isol} and Remark~\ref{rest} we obtain what we need.

Arguing as above, combining Lemma~\ref{local-direct-product} with
Theorem~\ref{add},
 we obtain
\begin{Cor}\hspace{-2mm}{\bf .}\label{coradd}
Let $L$ be a BM-structure on   connected  $(M^n,g)$. Suppose the
eigenvalue $\lambda_i$ has multiplicity $k$ at a typical point.
Suppose  there exists a point where the multiplicity of
$\lambda_i$ is greater than $k$. Then, there exists a point where
the multiplicity of $\lambda_i$ is precisely $k+1$.
\end{Cor}

Combining Lemma~\ref{local-direct-product} with
Corollary~\ref{typical}, we obtain

\begin{Cor}\hspace{-2mm}{\bf .}\label{9}
Let $L$ be a BM-structure on    connected  $(M^n,g)$.  Suppose the
eigenvalue $\lambda_i $ has multiplicity $k_i\ge 2$ at a typical point
and multiplicity $k_i+d$ at a point $p\in M^n$. Then, there exists
a point $q\in M^n$ in a small neighborhood of $p$  such that  the
eigenvalue $\lambda_i $ has multiplicity $k_i+d$ in $p$, and such
that $$ N_L(q)=\max_{x\in M^n}(N_L(x))-d.$$
\end{Cor}

We saw that under hypotheses of Theorems~\ref{isol},\ref{add},
the set of typical points is connected.
As it was shown in \cite{ERA}
(and will follow from Lemma~\ref{kol} from Section~\ref{instr}), in dimension 2 the set of typical points is connected as well. Combining these observations with   Lemma~\ref{local-direct-product}, we obtain

\begin{Cor}\hspace{-2mm}{\bf .}\label{010}
Let $L$ be a BM-structure on    connected  $(M^n,g)$ of dimension $n\ge 2$.
 Then, the set of typical points of $L$ is connected.
\end{Cor}

\subsection{ If $\phi_i$ is not isolated, $A_i$ has constant  sectional curvature. } \label{3.4}

In this section we assume that $L$ is a BM-structure on a
connected  complete Riemannian manifold
$(M^n,g)$. As usual, we denote by
$\lambda_1(x)\le...\le \lambda_n(x)$ the eigenvalues of $L$ at
$x\in M^n$.

\begin{Def}\hspace{-2mm}{\bf .} An eigenvalue $\lambda_i$ is called
{\bf isolated},  if,   for all  points $p_1,p_2\in M^n$
and for every $i,j\in \{1,...,n\}$, the equality
 $\lambda_i(p_1)=\lambda_j(p_1) $  implies
$\lambda_i(p_2)= \lambda_j(p_2) $.
\label{singul}
\end{Def}

As in Section~\ref{splitting},  at every point $p\in M^n$, we
denote by $V_i$ the eigenspace of $L$ with the eigenvalue
$\lambda_i(p)$. If $p$ is typical,  $V_i$ is a distribution near
$p$;
by Theorem~\ref{nijenhuis}, it is  integrable. We denote by
$M_i(p)$ the connected component (containing $p$)  of the
intersection  of the integral manifold with a small neighborhood
of $p$.

\begin{Th}\hspace{-2mm}{\bf .} \label{sing}
Suppose $\lambda_i$ is  a non-isolated eigenvalue such that its multiplicity at a typical point is greater than one.  Then, for every
typical point $p$, the restriction of $g$ to $M_i(p)$ has constant sectional
curvature.
\end{Th}

It could be easier to understand this Theorem using  the language
of Levi-Civita's Theorem~\ref{LC}: denote by
$\phi_1<\phi_2<...<\phi_m$  the different eigenvalues of $L $ at a
typical point. Theorem~\ref{sing}  says that if $\phi_i$ 
(of multiplicity $\ge 2$)  is
non-isolated, then $A_i$ from Levi-Civita's Theorem has constant sectional
curvature.

\noindent {\bf Proof of Theorem~\ref{sing}:}  Let $k_i>1$ be the  multiplicity  of $\lambda_i$ at a
typical point.
 Then, $\lambda_i$ is constant. Take a typical point
 $p$. We assume that $\lambda_i$ is not isolated; without loss of
 generality, we can suppose $\lambda_{i-1}(p_1)=\lambda_{i+k_i-1}(p_1)$ for some point
 $p_1$. By Corollary~\ref{coradd}, without loss of generality, we can
 assume  $\lambda_{i-1}(p_1)=\lambda_i(p_1)<\lambda_{i+k_i}(p_1)$.
 By Corollary~\ref{9},
  we can also assume that $N_L(p_1)=N_L(p)-1$, so that multiplicity of $\lambda_i(p_1)$ is precisely $k_i+1$.

Consider a geodesic segment  $\gamma:[0,1]\to M^n$ connecting $p_1$ and
$p$,\ $\gamma(0)=p$ and $\gamma(1)=p_1$.  Since it is sufficient
to prove Theorem~\ref{sing} at almost every typical point,
  without loss of
generality, we can assume that $p_1$ is the only nontypical point
of the geodesic segment $\gamma(\tau), \ \tau\in [0,1]$.
More precisely, take
$j\not\in i-2, i, ..., i+k_i-1$.  If  there
exists a point  $p_2\in M^n$ such that   
 $\lambda_j(p_2)<\lambda_{j+1}(p_2)$,
then, by assumptions,  $\lambda_j(p_1)<\lambda_{j+1}(p_1)$. Hence,
 by Lemma~\ref{technical},
for almost every  $\xi\in T_{\gamma(1)}M^n$, we have
$t_j(\gamma(1),\xi)\ne \max_{x\in M^n}(\lambda_j(x))$. Thus,
 almost all  geodesics starting
from $p_1$ do not contain  points where    $\lambda_j=\lambda_{j+1}$.
Finally, for almost every   $p\in M^n$, the geodesics connecting
 $p$ and $p_1$  contain  no  points where    $\lambda_j=\lambda_{j+1}$.

  Take the  point $q:=\gamma(1-\epsilon)$ of the segment,
where $\epsilon>0$ is small  enough. By Corollary~\ref{corisol},
the restriction of $g$ to $M_i(q)$ has constant sectional curvature.

Let us prove that  the geodesic segment $\gamma(\tau), \ \tau\in
[0,1-\epsilon]$  is  orthogonal to $M_i(\gamma(\tau))$ at every
point.
Consider the function $$\tilde I:TM^n\to \mathbb{R}; \
\tilde I(x,\xi):= \left(\frac{I_t(x,\xi)}{(\lambda_i-t)^{k_i-1}}
\right)_{|t=\lambda_i}. $$ Since the multiplicity of
$\lambda_i$
at every point  is at least $k_i$, by Lemma~\ref{technical},
 the function $\left(
\frac{I_t(x,\xi)}{(\lambda_i-t)^{k_i-1}} \right)$ is polynomial in
$t$ of degree $n-k_i$. Since $I_t$ is an integral,  for every fixed $t$, the function
$\left(\frac{I_t(x,\xi)}{(\lambda_i-t)^{k_i-1}}
\right)$   is an integral for the geodesic flow of $g$.
Thus, $\tilde I$ is
an  integral.

 At the tangent space to every point of geodesic
 $\gamma$, consider a coordinate system
 such that
 $L=\diag(\lambda_1,...,\lambda_n)$ and $g=\diag(1,...,1)$. In
 this coordinates,  $I_t(\xi)$ is given by (\ref{33}). Then,
 the integral $\tilde I(\xi)$ is the sum (we assume
 $\xi=(\xi_1,...,\xi_n)$)
 \begin{eqnarray} \label{33b}
\left( \sum_{\alpha=i}^{i+k_i-1}
\left(\xi_\alpha^2\prod_{\begin{array}{c}\beta=1
\\ \beta\ne
{i,i+1,...,i+k_i-1}\end{array}}^{n}(\lambda_\beta-\lambda_i)
\right)\right) \\ + \left( \sum_{\begin{array}{c}\alpha=1\\
\alpha\ne i,i+1,...,i+k_i-1\end{array}}^{n}
\left(\xi_\alpha^2\prod_{\begin{array}{c}\beta=1
\\ \beta\ne
{i+1,...,i+k_i-1}\\ \beta\ne \alpha
\end{array}}^{n}(\lambda_\beta-\lambda_i) \right)\right).
\label{33c}
 \end{eqnarray}

 Since the geodesic passes  through
the point where $\lambda_{i-1}=\lambda_i=...=\lambda_{i+k_i-1}, $
all products in the formulae above contain the factor
$\lambda_i-\lambda_i$, and, therefore, vanish,  so that $\tilde
I(\gamma(0),\dot\gamma(0))=0$. Since $\tilde I$ is an integral,
$\tilde I(\gamma(\tau),\dot\gamma(\tau))=0$ for every $\tau$. Let
us show that it implies that the geodesic is orthogonal to $M_i$
at every  typical point, in particular,  at points lying on the
segment  $\gamma(\tau), \ \tau\in [0,1[$.

Clearly, every term in the sum (\ref{33c}) contains the factor
$\lambda_i-\lambda_i$, and, therefore, vanishes. Then, the
integral $\tilde I$ is equal to (\ref{33b}).

At a typical point,   we have $$\lambda_1\le... \le
\lambda_{i-1}<\lambda_i=...=\lambda_{i+k_i-1}<\lambda_{k_i}\le...\le\lambda_n.$$
Then,
 the  coefficient  $$\prod_{\begin{array}{c}\beta=1
\\ \beta\ne
{i,i+1,...,i+k_i-1}\end{array}}^{n}(\lambda_\beta-\lambda_i) $$
is  nonzero. Then, all components
$\xi_\alpha$, $\alpha\in {i,...,i+k_i-1}$  vanish. Thus, $\gamma$ is
orthogonal to  $M_i$ at every typical point.

Finally, by Corollary~\ref{corlc}, the restriction of $g$  to
$M_i(p)$  is proportional to the restriction of $g$ to  $M_i(q)$
and, hence, has constant sectional  curvature.
 Theorem is
proved.

\subsection{If $g$ is $V(K)$-metric,  if $\phi_\aleph$ is not isolated, and
  if  the sectional  curvature of  $A_\aleph$
 is constant, then it is equal to $K_\aleph$. }  \label{3.5}
Within this section we assume that
 $L\ne \mathrm{const}\cdot \Id$ is  a BM-structure on a
 connected Riemannian manifold
$(M^n,g)$ of dimension $n\ge 3$. We denote by $m$ the number of different
eigenvalues of $L$ at a typical point.
  For every typical point $x\in M^n$, consider
 the Levi-Civita coordinates $(\bar x_1,..,\bar x_m)$
such that the metric  has the form (\ref{g}). We assume that there exists $i$ such that $k_i=1$.  Recall that the functions
$\phi_i$ are the eigenvalues of $L$:  in the Levi-Civita  coordinates,
$$
L=\diag(\lambda_1,...,\lambda_n)=\diag(\underbrace{\phi_1,...,\phi_1}_{k_1},....,
\underbrace{\phi_m,...,\phi_m}_{k_m}).$$

Consider  $\aleph\in \{1,...,m\}$ such that
$k_\aleph\ge 2$. We  put  $r:= k_1+...+k_{\aleph-1}$.
Let us assume that the eigenvalue $\phi_\aleph=\lambda_{r+1}$ 
is not isolated; that means that there exists a point $p_1\in M^n$
such that $\lambda_{r}(p_1)=\lambda_{r+1}$,
 or
 $\lambda_{r+k_{\aleph}+1}(p_1)=\lambda_{r+1}$.

Let us assume in addition that    in  a neighborhood
of every typical point,    the following holds:
\begin{enumerate}
\item The sectional curvature of   $A_\aleph$ is constant, \label{as1}
\item $g$ is a $V(K)$ metric. \label{as2}
 \end{enumerate}

As we saw in Section~\ref{3.4},  the  assumption \ref{as1}
 follows from the previous assumptions,
 if the  metric is complete.
As we saw in Section~\ref{Vk},  the  assumption \ref{as2} is
 automatically fulfilled, if the
 space of all  BM-structures is more than two-dimensional.

\begin{Th}\hspace{-2mm}{\bf .} \label{sing1} Under the  above  assumptions, the
sectional curvature of $A_\aleph$  is equal to $K_\aleph$.
\end{Th}

Recall that the definition of $K_\aleph$ is in the second statement of
Theorem~\ref{maximal}.

\noindent{\bf Proof of Theorem~\ref{sing1}:}
Let us denote by $\bar K_\aleph$ the sectional curvature of the metric $A_\aleph$.
By assumptions, it is constant
 in a neighborhood of every typical point.
Since the set of typical points is connected by Corollary\ref{010}, $\bar K_\aleph$ is independent
of the  typical point.
Similarly, since $K_\aleph$ is locally-constant by Theorem~\ref{maximal},
$K_\aleph$ is independent  of the typical  point.
Thus, it is sufficient to find a point where  $\bar K_\aleph=K_\aleph$.

Without loss of generality, we can suppose that there exists $p_1\in M^n$
such that $\lambda_{r}(p_1)=\lambda_{r+1}$.

 By Corollaries~\ref{9},\ref{010} without loss of generality we can assume that the multiplicity of $\lambda_{r+1}$ is $k_\aleph+1$
 in $p_1$,
 and that
  $N_L(p_1)=m-1$. Take a typical point $p$ in a small neighborhood of $p_1$.

Then, by Corollary~\ref{corisol}, the submanifold
 $M_{r+1}(p)$ is homeomorphic to the sphere. Since it is compact, there exists
a set of local coordinates charts on it such that there exist constants
$\textrm{const}>0$ and $\textrm{CONST}$ such that, in every
  chart $(x_\aleph^1,...,x_\aleph^{k_\aleph})$, for every
$\alpha, \beta\in \{1,..,k_\aleph\}$,   the entry
$(A_\aleph)_{\beta\beta}$ is greater than $\textrm{const}$ and the absolute
 value of the entry  $(A_\aleph)_{\alpha\beta}$ is less than   $\textrm{CONST}$,
  i.e.  $  A_\aleph(\frac{\partial }{\partial x_\aleph^\beta}, \frac{\partial }{\partial x_\aleph^\beta})> \textrm{const}$;
$ |A_\aleph(\frac{\partial }{\partial x_\aleph^\alpha}, \frac{\partial }{\partial x_\aleph^\beta})|< \textrm{CONST}$.

By shifting these local coordinates  along the vector fields
$\frac{\partial }{\partial x_i^j}$, where $i\ne \aleph$, for every typical
 point $p'$ in a neighborhood of $p_1$,
 we obtain coordinate charts on    $M_{r+1}(p')$ such that for every $\alpha,\beta$,
$(A_\aleph)_{\beta\beta}>\textrm{const}$,
 $|(A_\aleph)_{\alpha\beta}|<\textrm{CONST}$.

Let us calculate the projective Weyl tensor $W$ for $g$ in these local
coordinate charts.
We will be interested in the components (actually, in one component) of
$W$ corresponding to the coordinates $\bar x_\aleph$.
In what follows we reserve
the Greek letter $\alpha,\beta$ for the coordinates from $\bar x_\aleph$, so
that, for example,   $g_{\alpha\beta}$ will mean the component of the
metric staying on the intersection of  column number $r+\beta$   and  row
 number
$r+\alpha$.

As we will see below, the formulae  will include only the components
of $A_\aleph$. To simplify the notations, we will not write subindex $\aleph$ near $A_\aleph$, so
for example,   $g_{\alpha\beta}$ is equal to $P_\aleph \, A_{\alpha\beta}.$

Take $\alpha\ne \beta$.    Let calculate the component $W^\alpha_{\beta\beta\alpha}$. In order to do
it by formula~(\ref{pw}), we need to know $R^\alpha_{\beta\beta\alpha}$ and $R_{\beta\beta}$. It  is not easy
to calculate them: a straightforward way  is to
calculate  $R^\alpha_{\beta\beta\alpha}$ and $R_{\beta\beta}$ for the
metric (\ref{g}), then combine it  with   assumption \ref{as2}
(which could be written as a system of  partial differential equations) and
with  assumption \ref{as1} (which is a system of algebraic equations).
These was done
in \S 8 of \cite{solodovnikov1}, see formula (8.14) and what goes after it there. Rewriting the results of Solodovnikov in our notations, we obtain
$$
R^\alpha_{\beta\beta\alpha}=
\left(\bar K_\aleph  - (K_\aleph-K\, P_\aleph)\right) A_{\beta\beta}, \ \ R_{\beta\beta}=\left((k_\aleph-1)\, \bar K_\aleph + K\, (n-1)\, P_\aleph- (k_\aleph-1) \, K_\aleph  \right)\, A_{\beta\beta}.
$$
Substituting these expressions in the formula for projective  Weyl tensor
(\ref{pw}), we obtain
$$
W^\alpha_{\beta\beta\alpha}=(\bar K_\aleph-K_\aleph)\frac{n-k_\aleph}{n-1} A_{\beta\beta}.
$$
We see that, if $\bar K_\aleph\ne K_\aleph$, the component
$W^\alpha_{\beta\beta\alpha}$ is bounded from zero.

But if we consider a sequence of typical  points converging to $p_1$,
the component $W^\alpha_{\beta\beta\alpha}$ converge to zero. Indeed,
by definition $W^\alpha_{\beta\beta\alpha}=
W(dx_\aleph^\alpha,\frac{\partial }{\partial x_\aleph^\beta},\frac{\partial }{\partial x_\aleph^\beta},\frac{\partial }{\partial x_\aleph^\alpha})$, and the length of  $\frac{\partial }{\partial x_\aleph^\beta}$ goes to zero,
 if the  point goes to $p_1$.
Finally, $\bar K_\aleph= K_\aleph$.
Theorem is proved.

\section{Proof of Theorem~\ref{main}} \label{proof}

 If the dimension of ${\cal B}(M^n,g)$ is one,
 Theorem~\ref{main}
 is trivial: every projective transformation is
 a homothety.
  In Section~\ref{2proof},  we prove Theorem~\ref{main} under the additional assumption that the
   dimension of ${\cal B}(M^n,g)$ is two (Theorem~\ref{2dim}).
 In Section~\ref{3dim}, we prove Theorem~\ref{main} under the additional assumptions that
 the dimension of  $M^n$ and
  the dimension of ${\cal B}(M^n,g)$ are at least three
  (Theorem~\ref{4.2}). The last case, namely when the dimension
  of the manifold is two and
the dimension
  of  ${\cal B}$  is at least three requires
   essentially different methods and will be considered in Section~\ref{dim2}
  (Theorem~\ref{d2}).

\subsection{If the space ${\cal B}(M^n,g)$ has dimension two}
\label{1} \label{2proof}

Suppose $g$ and $\bar g$ are projectively equivalent. The
next lemma shows that the spaces ${\cal B}(M^n,g)$
and ${\cal
B}(M^n,\bar g)$ are canonically isomorphic:

\begin{Lemma}\hspace{-2mm}{\bf .}\label{10} Let $L$ be the BM-structure  (\ref{l})
constructed for the projectively equivalent  metrics $g$ and $\bar g$.
Suppose  $L_1$ is one more  BM-structure for $g$.
Then, $L^{-1}
L_1$ is a BM-structure for $\bar g$.
\end{Lemma}

\begin{Cor}\hspace{-2mm}{\bf .} \label{19}
If  ${\cal B}(M^n,g)$ is two-dimensional,
 every projective transformation
takes typical points  to typical points.
\end{Cor}
\noindent{\bf Proof of Lemma~\ref{10}:} It is sufficient to prove
the  statement locally. Let us fix a coordinate system and think
about tensors as about matrices. For every sufficiently big
constant $\alpha$,   the tensor $L+\alpha\cdot \Id$ is positive
defined. Then, by Theorem~\ref{th1}, $$ \bar{\bar g}:=\frac{1}{\det(L_1+\alpha\cdot
\Id)}\cdot g \cdot (L_1+\alpha\cdot \Id)^{-1} $$ is a Riemannian
metric projectively  equivalent to $g$. Then, it is projectively
equivalent to $\bar g$. Direct calculation of
  the tensor (\ref{l}) for the metrics $\bar g$, $\bar{\bar g}$ gives us that
$L^{-1} (L_1+\alpha\cdot \Id)$  is a  BM-structure for
$\bar g$. Since it is true for all big $\alpha$, and since
${\cal B}(M^n,\bar g)$ is a linear space, $L^{-1} L_1$ is a BM-structure for
$\bar g$. Lemma is proved.

\begin{Def}\hspace{-2mm}{\bf .}
A vector field $v$ is called {\bf projective}, if its flow consists of
projective transformations.
\end{Def}

A smooth one-parameter family of  projective transformations
$F_t:M^n\to M^n$  immediately gives us a projective vector field
$\left(\frac{d}{dt}F_t\right)_{|t=0}$. A projective vector field
gives us a one  parameter family of  projective transformations if
and only if it is complete.

\begin{Th}\hspace{-2mm}{\bf .} \label{2dim}
Suppose $(M^n, g)$ of dimension
 $n\ge 2$  is  complete and  connected.
Let $v$ be a complete projective vector field.
 Assume in addition that the dimension of  ${\cal B}(M^n, g)$
 is precisely  two.
Then, the flow of the vector field acts by affine transformations,
or the metric has constant positive sectional  curvature.
\end{Th}

\noindent{\bf Proof:}  Denote by $F_t$ the flow of the vector
field $v$. If it contains  not only affine transformations, there
exists $t_0\in \mathbb{R}$ such that the pull-back $\bar
g:=F_{t_0}^*(g)$ is projectively equivalent to $g$, and is not
affine equivalent to $g$.

Consider the BM-structure $L\in  {\cal B}(M^n, g)$ given by
(\ref{l}).  Take a  typical point $p$ such that $v$ does not
vanish at $p$ . Since $\bar g$ is not affine equivalent to $g$,
without loss of generality,  we can assume that  at least one
eigenvalue is not constant near $p$.   By Levi-Civita's Theorem
\ref{LC}, there exists a coordinate system $\bar x=(\bar
x_1,...,\bar x_m)$ in a neighborhood of $p$ such that the tensor
$L$ and  the metrics $g,\bar g$ have the form
(\ref{diagonal},\ref{g},\ref{rem2}), respectively.
 In particular, all these objects are block-diagonal
with the parameters of the blocks $(k_1,...,k_m)$. Note that the
nonconstant eigenvalue of $L$  has multiplicity one.

Let  the
vector $v$ be  $(\bar v_1,...,\bar v_m)$ in this coordinate system.
Then, by Theorem~\ref{projective} and Lemma~\ref{10},
  and using that the space ${\cal B}(M^n, g)$ is two-dimensional,
we obtain that the Lie derivatives ${\cal L}_vg$ and  ${\cal
L}_v\bar g$ are block-diagonal as well.  In the coordinate system
$\bar x=(\bar x_1,...,\bar x_m)$,  the metric $\bar g$ is given by
the matrix $\frac{1}{\det (L)}gL^{-1}$. Then, the  Lie derivative
${\cal L}_v\bar g$
 is $$\left({\cal L}_v\frac{1}{\det
(L)}\right)gL^{-1}+\frac{1}{\det (L)} \left({\cal
L}_vg\right)L^{-1}+\frac{1}{\det (L)}g \left({\cal L}_v
L^{-1}\right). $$ We see that the first two terms of the sum are
block-diagonal; then, the third term must be block-diagonal as
well so that ${\cal L}_v L$ is block-diagonal. Let us calculate
the element of ${\cal L}_v L$ which is on the intersection of
$x_i^\alpha$-row  and $x_j^\beta$-column. If $i\ne j$, it is equal
to $$ \pm(\phi_i-\phi_j)\frac{\partial v_i^\alpha}{\partial
x_i^\beta}. $$ Thus, $$ \frac{\partial v_i^\alpha}{\partial
x_i^\beta}\equiv 0. $$
 Finally,  the block $\bar v_i$ of the projective vector field
depends on the variables $\bar x_i$ only.

Using this, let us calculate  the tensor $g^{-1}{\cal L}_vg$
 from Theorem~\ref{projective}
for the metric $g$.   We see that (in matrix notations), the
matrix of $g$  is $AP$, where $$
 P:=\diag(\underbrace{P_1,...,P_1}_{k_1},....,
\underbrace{P_m,...,P_m}_{k_m}) \ \   \textrm{and} \ \
A:=\textrm{block-diagonal}(A_1,...,A_m).
$$
We understand $A$ as a (0,2)-tensor and  $P$ as a (1,1)-tensor. Then,
$$
 g^{-1}{\cal L}_vg=P^{-1}A^{-1}({\cal L}_vA)P+P^{-1}{\cal L}_vP.
$$ Since the entries of $A_i$ and of $\bar v_i$  depend on the
coordinates $\bar x_i$  only, ${\cal L}_v A$ is block-diagonal so
that the first term on the right hand side is equal to
$A^{-1}{\cal L}_vA$. Thus,
\begin{equation} \label{Lg}
 g^{-1}{\cal L}_vg=A^{-1}{\cal L}_vA+P^{-1}{\cal L}_vP.
\end{equation}
 By direct calculations, we obtain:
\begin{eqnarray*}
P^{-1}{\cal L}_vP &=& \diag\left( \underbrace{\sum_{i=2,...,m} \frac{\phi_1'v_1^1}{\phi_i-\phi_1}+ \frac{\phi_i'v_i^1}{\phi_1-\phi_i},...,\sum_{i=2,...,m} \frac{\phi_1'v_1^1}{\phi_i-\phi_1}+ \frac{\phi_i'v_i^1}{\phi_1-\phi_i}}_{k_1},\right. \\
&&..., \left.
\underbrace{\sum_{i=1,...,m-1} \frac{\phi_m'v_m^1}{\phi_i-\phi_m}+ \frac{\phi_i'v_i^1}{\phi_m-\phi_i},...,\sum_{i=1,...,m-1} \frac{\phi_m'v_m^1}{\phi_i-\phi_m}+ \frac{\phi_i'v_i^1}{\phi_m-\phi_i}}_{k_m}
   \right) .
\end{eqnarray*}

By Theorem~\ref{projective}, since the space ${\cal B}(M^n,g)$ is
two-dimensional,     $g^{-1}{\cal L}_vg$  equals
 $$
 \alpha L+ \beta\cdot \Id + \textrm{trace}_{\alpha L+ \beta\cdot \Id}\cdot \Id,
$$
where $\alpha,\beta \in \mathbb{R}$.
Then,  the first entry of  the block number $j$ gives  us the
following equation:
 \begin{equation} \label{block}
  a_j(\bar x_j)+\sum_{\begin{array}{lcr}i&=&1,...,m\\i&\ne& j\end{array}}
\frac{\phi_j'v_j^1}{\phi_i-\phi_j}+
\frac{\phi_i'v_i^1}{\phi_j-\phi_i}=
\alpha\phi_j+\beta+\textrm{trace}_{\alpha L+ \beta\cdot \Id},
\end{equation}
 where $a_j(\bar x_j)$ collects the terms coming from
$A^{-1}{\cal L}_vA$,  and, hence,  depends on the variables $\bar
x_j$ only.

Our next goals is to show that

\begin{itemize}
\item Only one $\phi_j$ can be nonconstant. The behavior of  nonconstant  $\phi_j$
on the orbit passing through $p$ is   given by  $$ \phi_j(F_t(p))
:=  \frac{b}{2\alpha }+  \frac{\sqrt{D}}{\alpha }\tanh(\sqrt{D}(t+
d)), $$ where $D:= b^2/4+\alpha c$, where $b, \ c $ and $d$ are
certain (universal along the orbit) constants.

\item The constant eigenvalues $\phi_s$ are  roots of  the polynomial
$$-\alpha\phi^2 +b\phi +c.$$ In particular, $m$ is at most 3.
\end{itemize}

Take $s\ne j\in 1,...,m$. We see that the terms  in (\ref{block})
depending on the variables $\bar x_s$ are $$
\Phi_j:=\frac{\phi_j'v_j^1}{\phi_s-\phi_j}+
\frac{\phi_s'v_s^1}{\phi_j-\phi_s}- k_s\alpha\phi_s. $$ Thus, $
\Phi_j$ depends on the variables $\bar x_j$ only. Similarly, $$
\Phi_s:=\frac{\phi_s'v_s^1}{\phi_j-\phi_s}+
\frac{\phi_j'v_j^1}{\phi_s-\phi_j}- k_j\alpha\phi_j $$ depends on
the variables $\bar x_s$ only. Using that  $ \Phi_j- \Phi_s$ is
equal to
$
\alpha(k_j\phi_j-k_s\phi_s)$, we obtain that (for an appropriate
constant $B$) $\Phi_j$ must be equal to $\alpha k_j \phi_j+ B$.
Thus,
\begin{equation} \label{nc}
\frac{\phi_j'v_j^1}{\phi_s-\phi_j}+
\frac{\phi_s'v_s^1}{\phi_j-\phi_s}- k_s\alpha\phi_s=\alpha k_j
\phi_j+ B.
\end{equation}

Now let us prove that at least one of the functions  $\phi_s$ and
$\phi_j$ is constant near $p$. Otherwise,  $k_s=k_j=1$, and
Equation (\ref{nc}) is equivalent to $$ \phi_j'v_j^1+\alpha
\phi_j^2-B \phi_j = \phi_s'v_s^1+\alpha \phi_s^2-B \phi_s. $$ We
see that the terms on the left-hand side depend on the variable
$x_j^1$ only, the terms on the right-hand side depend on the
variable $x_s^1$ only. Then, there exists a constant $c$ such that
\begin{eqnarray*}
\phi_j'v_j^1 & = & -\alpha \phi_j^2+b \phi_j+c \\ \phi_s'v_s^1 & =
& -\alpha \phi_s^2+b \phi_s+c,
\end{eqnarray*}
where $b:=B$.
 Since $\phi_j$ and $v_j^1$ depend on the variable
$x_j^1$ only, $\phi_j'v_j^1$  at the point $F_t(p)$ is equal to
$\dot \phi_j:=\frac{d}{dt}\phi_j(F_t(p))$, so that
$\phi_s(F_t(p))$ and $\phi_j(F_t(p))$ are  solutions of the
following differential equation:
\begin{equation}\label{de}
 \dot \phi  =  -\alpha \phi^2+b \phi+c.
 \end{equation}
By Corollary~\ref{19}, all  points $F_t(p)$ are typical. Then,
$\phi_s$ and $\phi_j$ are solutions of  Equation (\ref{de}) at
every point of the orbit passing through $p$, and the constants
$\alpha, b, c$ are universal along the orbit.
 The
equation can be solved (we assume  $\alpha\ne 0$). The
  constant solutions are
 $\frac{b}{2\alpha }\pm \frac{\sqrt{{b^2}/{4}+{\alpha
c}}}{\alpha}$, and  the  nonconstant solutions are:
\begin{enumerate}

\item For  $D:={b^2/4+\alpha c}<0 $, every nonconstant  solution is the  function

  $\frac{b}{2\alpha}+ \frac{\sqrt{-D}}{\alpha}\tan(\sqrt{-D}(t+ d_1)).$ \label{f1}

\item
For  $D:={b^2}/{4}+{\alpha c}>0 $, every nonconstant  solution is
one of the  functions
\begin{enumerate}
\item $\frac{b}{2\alpha }+  \frac{\sqrt{D}}{\alpha }\tanh(\sqrt{D}(t+ d)), $ \label{r}
\item $\frac{b}{2\alpha }+ \frac{\sqrt{D}}{\alpha }\coth(\sqrt{D}(t+ d_3))$
\label{f2}
\end{enumerate}

\item
For  $D:={b^2}/{4}+{\alpha c}=0 $, every nonconstant  solution is
the  function

 $\frac{b}{2\alpha }+\frac{1}{\alpha (t+ d_4)}. $ \label{f3}

\end{enumerate}

 The solutions  (\ref{f1},\ref{f2},\ref{f3}) explode
in finite  time. This gives us a contradiction: the metrics $g$
and $\bar g$ are smooth and, hence,    the eigenvalues of $L$ are
finite on every compact set.

If the functions $\phi_s(F_t(p))$  and $ \phi_j(F_t(p))$ have the
form (\ref{r}), then there exist points  $p_1,p_2,q_1,q_2$ of the
orbit passing through $p$, such that $\phi_s(p_1)<\phi_j(p_2)$,
$\phi_s(q_1)>\phi_j(q_2)$. This gives a contradiction with
Corollary~\ref{ordered12}.

Thus, only one eigenvalue of $L$  can be nonconstant in  a
neighborhood of $p$. Let $\phi_j$ be nonconstant near $p$.

Now let us show that a constant eigenvalue is  a  root of the
polynomial $-\alpha \phi^2+b\phi + c$; in particular, there are no
more than two different constant eigenvalues.

Suppose $\phi_s$ is constant. Then,  the  derivative $\phi_s'$
vanishes and (\ref{nc}) reads: $$ \phi_j'v_j^1=-\alpha
\phi_j^2+(\alpha \phi_s- \alpha k_s\phi_s-
B)\phi_j+\alpha\phi_s^2k_s+\phi_s B. $$ Denoting $\alpha \phi_s-
\alpha k_s\phi_s- B$ by $b$ and $\alpha\phi_s^2k_s$ by $c$, and
arguing as above, we see that $\phi_j$ is a solution of
(\ref{de}).  Hence,  the  behavior of  $\phi_j$ on the orbit
(passing through $p$) is given by (\ref{r}). Clearly, $\phi_s$ is
a root of ${-\alpha \phi^2+b \phi+c}$.

Thus, near $p$  only the following three case are possible:

Case 1: $m=3$. The eigenvalues $\phi_1,\ \phi_3$ are constant; the
eigenvalue $\phi_2$ is not constant.

Case 2a: $m=2$. The eigenvalues $\phi_1$ is  constant; the
eigenvalue $\phi_2$ is not constant.

Case 2b:  $m=2$. The eigenvalues $\phi_2$ is  constant; the
eigenvalue $\phi_1$ is not constant. \vspace{2ex}

In all three cases, one can prove that the metric has constant sectional
curvature.  We will carefully consider the most complicated case 1
and   sketch the proof for case 2a. The proof for case 2b is
similar to  the proof for case 2a.

Suppose $m=3$,  the eigenvalues  $\phi_1,\ \phi_3$ are constant and
the eigenvalue $\phi_2$ is not constant in a neighborhood of $p$.
Without loss of generality, we can assume $A_2(dx_2^1,dx_2^1)=(dx_2^1)^2$.
  Then, $\alpha L+ \beta\cdot \Id$ is

$$
\diag\left(\underbrace{\alpha\phi_1+\beta,...,\alpha\phi_1+\beta}_{k_1},
\alpha\phi_2+\beta,\underbrace{\alpha\phi_3+\beta,...,\alpha\phi_3+\beta}_{k_3}\right).
$$
Since $g^{-1}{\cal L}_vg$ is equal to  $\alpha L+ \beta\cdot \Id+
\mathrm{trace}_{\alpha L+ \beta\cdot \Id}\cdot \Id$,
  Equation~(\ref{Lg}) gives us the following system:
\begin{equation}\label{sys1}
\hspace{-5mm}\left\{ \hspace{-3mm} \begin{array}{cl}
(\alpha\phi_1+\beta)(k_1+1)+\alpha\phi_2+\beta+(\alpha\phi_3+\beta)k_3&=a_1(\bar
x_1)+\frac{\phi_2'v_2^1}{\phi_1-\phi_2}\\
(\alpha\phi_1+\beta)k_1+2(\alpha\phi_2+\beta)+(\alpha\phi_3+\beta)k_3&=-2
\frac{\partial v_2^1}{\partial x_2^1} +
\frac{\phi_2'v_2^1}{\phi_1-\phi_2}
+\frac{\phi_2'v_2^1}{\phi_3-\phi_2}\\
(\alpha\phi_1+\beta)k_1+\alpha\phi_2+\beta+(\alpha\phi_3+\beta)(k_3+1)&=a_3(\bar
x_3)+\frac{\phi_2'v_2^1}{\phi_3-\phi_2}.
\end{array}\right.
\end{equation}
Here $a_{1}, \ a_3$ collect the terms coming from $A_1^{-1}{\cal
L}_vA_1$ and $A_3^{-1}{\cal L}_vA_3$, respectively.   Using
$\phi_2'v_2^1=\dot \phi_2 =
-\alpha(\phi_1-\phi_2)(\phi_3-\phi_2)$, we obtain

\begin{equation}\label{sys2}
\left\{ \begin{array}{ccc}
(\alpha\phi_1+\beta)(k_1+1)+(\alpha\phi_3+\beta)(k_3+1)&=&a_1(\bar
x_1)\\ -\frac{a_1(\bar x_1)}{2}&=& \frac{\partial v_2^1}{\partial
x_2^1}
\\
(\alpha\phi_1+\beta)(k_1+1)+(\alpha\phi_3+\beta)(k_3+1)&=&a_3(\bar
x_3).
\end{array}\right.
\end{equation}

We see that $a_1=a_3=const$; we denote this constant by $a$. Let
us prove that $a=0$. We assume that $a\ne 0$ and will find a
contradiction.

Consider  Equation (\ref{Lg}).
 We see that the first block of the left-hand side and the first block of the
 second term at the right-hand side are proportional to  $\diag(\underbrace{1,...,1}_{k_1})$. Then, $A_1^{-1}{\cal
 L}_vA_1$ is proportional to identity. The coefficient of
 proportionality is clearly $a$. Then, $A_1^{-1}{\cal
 L}_vA_1= a \diag(\underbrace{1,...,1}_{k_1})$, and ${\cal
 L}_vA_1=aA_1$, so that the flow of the vector field $(\bar v_1,
 0, \underbrace{0,...,0}_{k_3})$ acts by homotheties on the
 restriction of $g$ to the coordinate plaque of the coordinates
 $\bar x_1$. Note that this vector field is the orthogonal
 projection of $v$ to the coordinate plaque.

 Similarly, the vector field $(
 \underbrace{0,...,0}_{k_1},0,\bar v_3)$ acts by homotheties  (with the same coefficient of stretching)
 on the
 restriction of $g$ to the coordinate plaque of the coordinates
 $\bar x_3$.

 Without loss of generality, we can suppose that $\phi_2(p)$ is a regular value of $\phi_2$.
 Denote by $\check M_2$  the connected component of the set $\{q\in M^n:
 \ \phi_2(q)=\phi_2(p)\}$ containing $p$. By construction,   $\check M_2$ is a
 submanifold of codimension $1$, and the derivative of $\phi_2$
 vanishes at no points of $\check M_2$. Then, at   every point of $\check
 M_2$,
 the flow of
 the orthogonal projection of $v$ to $\check M_2$ acts by homotheties.
 Since $M^n$ is complete, $\check M_2$ is complete as well. Then,
 $\check M_2$ with the restriction of the  metric $g$ is isometric to the
 standard Euclidean space $(\mathbb{R}^{n-1}, g_{euclidean})$, and
 there exists  precisely one point where the orthogonal projection
 of $v$ vanishes, see \cite{lichnerowicz} for details. Without
 loss of generality, we can think that $p$ is the point where the orthogonal projection
 of $v$ vanishes.

 Moreover, at every point of $\check M_2$, consider the
 eigenspaces of $L$ corresponding to $\phi_1 $ and $\phi_3$. By
 Theorem~\ref{nijenhuis}, they are tangent to $\check M_2$; by
 Corollary~\ref{ordered12}, every point of $\check M_2$ is
 typical so that the eigenspaces corresponding to $\phi_1 $ and
 $\phi_3$ give us two distributions. This distributions are
 integrable by Theorem~\ref{nijenhuis}.  We denote by $M_1(p)$ and
 $M_3(p)$ their integral manifolds passing through $p$. Locally,
 in Levi-Civita's coordinates from Theorem~\ref{LC}, the manifold
 $M_1$ is the  coordinate plaque of coordinates $\bar x_1$ and
 $M_3$ is the  coordinate plaque of coordinates $\bar x_3$. Then,
 $M_1(p)$ and $M_3(p)$ are invariant with respect to the orthogonal projection of $v$ to
 $\check M_2$.

 Consider the orbit of the projective action of $(\mathbb{R},+)$
 containing the point $p$. Since at the point $p$  the vector
  $v$ has the form
  $(\underbrace{0,...,0}_{k_1},v_2^1,\underbrace{0,...,0}_{k_3})$,
  and since the components $\bar v_1, \ \bar v_3$ do not depend
  on the coordinate $\bar x_2$, at every point of the orbit $v$ is
  an eigenvector of $L$ with the eigenvalue $\phi_2$.

  Let us show that the length of the orbit is finite at least  in one
  direction. Indeed, the second equation in (\ref{sys2}) implies
$$ \frac{d}{d t}\log\left( v_2^1(F_t(p))\right) = -\frac{a}{2}.  $$ Its solution is
$v_2^1(F_t(p))= \textrm{Const}  \exp\left(-\frac{a}{2}t\right)$. Then,
  the length between points $t_1<t_2$ of the orbit
  is equal to $$\int_{t_1}^{t_2}\sqrt{g(v,v)}\ dt = {\textrm{Const}}
  \int_{t_1}^{t_2}\sqrt{(\phi_2-\phi_1)(\phi_3-\phi_2)}\exp\left(-\frac{a}{2}t\right)\
  dt.
$$ Since $(\phi_2-\phi_1)(\phi_3-\phi_2)$ is bounded, the length
of the orbit  is finite at least in one direction.

Then, since the manifold is complete, the  closure of the orbit
contains   a point $q\in M^n$ such that either $ \phi_2(q)=\phi_1$
or $\phi_2(q)=\phi_3$. Without loss of generality, we can assume
that $ \phi_2(q)=\phi_1$.

 Then, without loss of generality, we can  assume that $p$ is close
 enough to $q$, so that, by Corollary~\ref{corisol}, \ $M_1(p)$ is
 homeomorphic to the sphere.
     We got a contradiction with the fact that  $M_1(p)$ admits a vector field whose flow acts by homotheties, see \cite{lichnerowicz}.
     Finally, $a=0$.

Since $a=0$, the second equation in (\ref{sys2}) implies $$
\frac{d}{d t}\log \left(v_2^1(t)\right) = 0 , $$ so that locally
\begin{equation} \label{2222}
x_2^1=
\textrm{Const}\ t.\end{equation}
 Then, the  adjusted metric has the form $$
(1-\tanh(y_2+d))\ (dy_1)^2 + C(1-\tanh(y_2+d))(1+\tanh(y_2+d))(dy_2)^2
+ (1+\tanh(y_2+d))(dy_3)^2 $$
 in a certain coordinate system and for
certain constants $C,\ d$. By direct calculation,  we see that the sectional
curvature  of the adjusted metric is positive constant. If the dimension of
the manifold is three, it implies that the sectional curvature of $g$ is constant.

If the dimension of the manifold is greater than three, in view of
Theorems~\ref{sing},\ref{sing1},\ref{maximal}, it is sufficient to   show that there exist points $q_1$, $q_3$
 such that
$\phi_2(q_1)=\phi_1$ and $\phi_2(q_3)=\phi_3$. Take the geodesic
$\gamma$ such that $\gamma(0)=p$ and
$\dot\gamma(0)=(\underbrace{0,...,0}_{k_1},v_2^1,\underbrace{0,...,0}_{k_3})$.

 Let us show that at every typical point of the geodesic, in Levi-Civita's coordinates,
 the $\bar x_1$- and $\bar x_3$-components  of the geodesic
 vanishes.

   Consider the functions
$$I':=\left(\frac{I_t}{(\phi_1-t)^{k_1-1}}\right)_{|t=\phi_1}:TM^n\to
\mathbb{R} ,$$
$$I'':=\left(\frac{I_t}{(\phi_3-t)^{k_3-1}}\right)_{|t=\phi_3}:TM^n\to
\mathbb{R}. $$ They are integrals of the geodesic flow. Since
$I'(\gamma(0),\dot\gamma(0))=I''(\gamma(0),\dot\gamma(0))=0$, at
every point $\tau$ of the geodesic we have  $
I'(\gamma(\tau),\dot\gamma(\tau))=I''(\gamma(\tau),\dot\gamma(\tau))=0.$
Then, at every typical point of the geodesic, the $\bar x_1$- and
$\bar x_3$- components of  the  velocity vector vanishes. Consider
the maximal (open)  segment of this geodesic containing $p$ and
containing only typical points. Let us show that this segment has
finite length; that $\phi_1=\phi_2$ at one
 end of the segment and $\phi_3=\phi_2$ at the other end of the
 segment.

Using (\ref{2222}), we obtain that
 $v_2^1=\textrm{Const}\cdot x^1_2$  near every    every point of the segment.
 Then, we can globally parameterize the coordinate $x_2^1$ near the  points
of the  geodesic segment such that the constant $\textrm{Const}$ is
universal along the segment. Then, the length of the segment  is
given by (we denote by $v_2$ the projection of $v$ to the segment)
$$ \int_{-\infty}^{+\infty}\sqrt{g(v_2,v_2)}\ dt = \textrm{Const}
  \int_{t_1}^{t_2}\sqrt{(\phi_2-\phi_1)(\phi_3-\phi_2)} \ t \
  dt.
$$ Since $(\phi_2-\phi_1)(\phi_3-\phi_2)$ decrees exponentially
for $t\longrightarrow  \pm \infty$, the length of the segment is
finite. Clearly, the limit of $\phi_2$ is $\phi_1$ in one
direction and $\phi_3$ in the other direction. Finally,
 there exists a   point where $\phi_1=\phi_2$ and a point where
$\phi_2=\phi_3$.

 Then, all eigenvalues of $L$ are not
 isolated.   Then, by Theorems~\ref{sing}, \ref{sing1},
 every $A_i$  has constant sectional  curvature $K_i$, and,
 by Theorem~\ref{maximal},
 $g$ has
 constant sectional curvature. By Corollary~\ref{bel}, the sectional
curvature of $g$ is positive.
Theorem~\ref{2dim} is proved under the assumptions of case 1.

The proof for cases 2a,2b is similar; we will sketch the proof for
case  2a: First of all, under the assumptions of case 2a,
 $\alpha L+ \beta\cdot \Id$ is $$
\diag\left(\underbrace{\alpha\phi_1+\beta,...,\alpha\phi_1+\beta}_{k_1},
\alpha\phi_2+\beta\right). $$ Since $g^{-1}{\cal L}_vg$ is equal
to  $\alpha L+ \beta\cdot \Id+\textrm{trace}_{\alpha L+\beta\cdot \Id}$,  Equation~(\ref{Lg}) gives us the
following system:
\begin{equation}\label{sys1a}
\left\{ \begin{array}{ccc}
(\alpha\phi_1+\beta)(k_1+1)+\alpha\phi_2+\beta&=&a_1(\bar
x_1)+\frac{\phi_2'v_2^1}{\phi_1-\phi_2}\\
(\alpha\phi_1+\beta)k_1+2(\alpha\phi_2+\beta)&=&-2 \frac{\partial
v_2^1}{\partial x_2^1} + \frac{\phi_2'v_2^1}{\phi_1-\phi_2}
\end{array}\right.
\end{equation}
As we have proven, $\phi_1$ is a root  of  ${-\alpha \phi^2+b
\phi+c}$.  We denote by $\phi_3$ the second root of ${-\alpha
\phi^2+b \phi+c}$. Arguing as in case 1, we have

\begin{equation}\label{sys2a}
\left\{ \begin{array}{ccc}
(\alpha\phi_1+\beta)(k_1+1)+(\alpha\phi_3+\beta)&=&a_1(\bar x_1)\\
\alpha(\phi_2-\phi_1)+{a_1(\bar x_1)}{}&=& -2\frac{\partial
v_2^1}{\partial x_2^1}
\end{array}\right.
\end{equation}
Then, $a_1=Const$. Arguing as in case 1, one can prove that
$a_1=0$. Then,  the second equation of (\ref{sys2a}) implies $$
\frac{d}{dt}  \log \left(v^1_2\right)= \alpha(1+\tanh(\sqrt{b^2/4+\alpha c}\
t+d)).$$

Thus, $v^1_2(t)=\frac{\exp(\sqrt{b^2/4+\alpha c}\
t/2+d/2)}{\sqrt{\cosh(\sqrt{b^2/4+\alpha c}\ t+d)}}$. Then, the
adjusted  metric  is  (locally) proportional to
  $$
  (1-\tanh(\sqrt{b^2/4+\alpha c}\ y_2+d))\left(dy_1^2+
  \frac{\exp(\sqrt{b^2/4+\alpha
c}\ y_2+d)}{\cosh(\sqrt{b^2/4+\alpha c}\ y_2+d)} \ dy_2^2\right)
 $$
 and, therefore, has
   constant  curvature. If $n=2$, it implies that $g$ has constant curvature.
If $n\ge 3$,   similar to the proof for  case 1,
 we can show that there exists a point where $\phi_1=\phi_2$.
Then, by Theorems~\ref{sing},\ref{sing1},\ref{maximal},
 the sectional curvature of $g$ is constant. By Corollary~\ref{bel},
it is positive.
Theorem~\ref{2dim} is proved.

\subsection{Proof if  $\dim(M^n)\ge 3$; $\dim({\cal B}(M^n,g))\ge
3$. } \label{3dim} Assume $\dim({\cal B}(M^n,g))\ge 3$, where
$(M^n,g)$ is a connected complete Riemannian metric of dimension
$n\ge 3$. Instead of proving Theorem~\ref{main} under these
assumptions, we will prove the following stronger
\begin{Th}\hspace{-2mm}{\bf .}  \label{4.2} { \it
Let $(M^n,g)$ be  a connected complete Riemannian manifold of
dimension $n\ge 3$. Suppose $\dim({\cal B}(M^n,g))\ge 3$.

Then, if a complete Riemannian metric $\bar g$ is projectively
equivalent to $g$, then  $g$ has positive constant sectional curvature or
$\bar g$ is affine equivalent to $g$. }
\end{Th}

\noindent {\bf Proof: }   Denote by $L$ the BM-structure from
Theorem~\ref{th1}.
In view of
Remark~\ref{eisen},  without loss of generality,
we can assume that at least one eigenvalue of $L$ is not constant.
\weg{Then, by Theorem~\ref{solodovnikov}, for every
typical point,  the sectional curvature of
 the adjusted metric is constant.}

Denote by $m$ the number of different eigenvalues of $L$ in a
typical point. The number $m$ does not depend on the typical
point. If $m=n$, Theorem~\ref{4.2} follows from
  Fubini's Theorem~\ref{fubini} and Corollary~\ref{bel}.

  Thus, we can assume $m<n$. Denote
by $m_0$ the number of simple eigenvalues of $L$ at a typical
point. By Corollary~\ref{ordered12}, the number $m_0$ does not
depend on the typical point. Then, by Levi-Civita's
Theorem~\ref{LC}, the metric $g$ has the following warped
decomposition near every typical point $p$: \begin{equation}
\label{wrap}
 g=g_0+\left|\prod_{\begin{array}{cc}i=1\\ i\ne m_0+1 \end{array}}^{m}
 (\phi_{m_0+1}-\phi_i)\right|g_{m_0+1}+...
 +\left|\prod_{\begin{array}{cc}i=1\\ i\ne m \end{array}}^{m}(\phi_{m}-\phi_i)\right|g_{m}.
\end{equation}
Here the coordinates  are $(\bar y_0,...,\bar y_m)$, where $\bar
y_0=(y_0^1,...,y_0^{m_0})$  and for $i>1$ \ \ $\bar
y_i=(y_i^1,...,y_i^{k_i})$. For $i\ge 0$, every  metric $g_{m_0+i}$
depends on the coordinates $\bar y_i$ only. Every  function
$\phi_i$
 depends on $y_0^i$ for  $i\le m_0$ and is  constant  for $i> m_0$.

 Let us explain the relation between Theorem~\ref{LC} and
the formula above.   The term  $g_0$ collects
all one-dimensional
terms of (\ref{g}).  The coordinates $\bar y_0$ collect all
one-dimensional $\bar x_i$ from (\ref{g}).  For $i>m_0$,  the
coordinate $\bar y_i$  is   one of the coordinates  $\bar x_j$
with $k_j>1$. Every metric $g_{m_0+i}$ for $i>1$ came from one of
the  multidimensional terms of (\ref{g}) and is proportional  to
the corresponding  $A_j$. The functions $\phi_i$ are eigenvalues
of $L$; they must not be ordered anymore: the indexing can be
different from (\ref{diagonal}). Note that,  by
Corollary~\ref{ordered12}, this re-indexing can be done
simultaneously in all  typical points.

Since the dimension of the  space ${\cal B}(M^n,g)$ is greater
than two, by Theorem~\ref{solodovnikov}, $g$ is a $V(K)$ metric
near every typical point. By Corollary~\ref{010}, the set of the
typical points is connected, so that the constant $K$ is
independent of the typical point.

 According to Definition~\ref{singul},
 a multiple eigenvalue
$\phi_i$ of $L$ is  {\bf isolated}, if there exists no nonconstant
eigenvalue $\phi_j$ such that $\phi_j(q)=\phi_i$ at some point
$q\in M^n$. If every multiple eigenvalue of $L$ is non-isolated,
then, by Theorems~\ref{sing},\ref{sing1},\ref{maximal},
  $g$ has constant
  sectional curvature. By Corollary~\ref{bel}, the sectional
curvature is positive.

 Thus, we
can assume that there exist isolated eigenvalues. Without loss of
generality, we can assume that (at every typical point) the
re-indexing of $\phi_i$  is made in such a way that the first
multiple eigenvalues $\phi_{m_0+1},...,\phi_{m_1}$ are
non-isolated and the last multiple eigenvalues
$\phi_{m_1+1},...,\phi_{m}$ are isolated. By assumption, $m_1< m$.

We will prove that in this case all eigenvalues of $L$ are
constant. By Remark~\ref{eisen}, it implies that the metrics $g$,
$\bar g$ are affine equivalent.

Let us show that $K$ is nonpositive.
 We  suppose that it is positive and will find  a
contradiction.

At every point $q$ of $M^n$, denote by $V_0\subset T_qM^n$  the
direct product of the eigenspaces of $L$ corresponding to the
eigenvalues $\phi_1,...,\phi_{m_1}$. Since the  eigenvalues
$\phi_{m_1+1},...,\phi_m$ are isolated by the assumptions, the
dimension of $V_0$ is constant, and $V_0$ is a distribution. By
Theorem~\ref{nijenhuis}, $V_0$ is integrable. Take a typical point
$p\in M^n$ and  denote by $M_0$ the integral manifold of the
distribution containing this point. Since $M_0$ is totally
geodesic,  the restriction $g_{|M_0}$ of the metric $g$  to $M_0$
is complete. By Theorems~\ref{sing},\ref{sing1},\ref{maximal}, the
metric $g_{|M_0}$ has constant sectional curvature $K$, or $M_0$ is one-dimensional.

Consider the direct product $M_0\times \mathbb{R}^{m-m_1}$ with
the metric \begin{equation} \label{g0}
g_{|M_0}+\left|\prod_{\begin{array}{c} i=1 \\ i\ne m_1+1
\end{array} }^{m}
 (\phi_{m_1+1}-\phi_i)\right|dt_{m_1+1}^2 + ... +\left|\prod_{\begin{array}{c}i=1 \\ i\ne m \end{array}}^{m}
 (\phi_{m}-\phi_i)\right|dt_m^2 , \end{equation} where  $(t_{m_1+1},...,t_m)$ are  the standard
 coordinates
on $\mathbb{R}^{m-m_1}$. Since the eigenvalues
$\phi_{m_1+1},...,\phi_m$ are isolated, (\ref{g0})  is a
well-defined Riemannian metric. Since $g_{|M_0}$ is complete, the
metric  (\ref{g0}) is complete. \weg{By definition, the metric is
the adjusted metric for the warped decomposition (\ref{wrap}).}
Since the sectional curvature of the adjusted metric is $K$, 
 and since the sectional curvature of $g_{|M_0}$ is $K$ (or $M_0$ has 
dimension $1$), the
 sectional curvature of (\ref{g0}) is $K$ as well.
  If $K> 0$, then  the
product $M_0\times \mathbb{R}^{m-m_1}$ must be compact, which
contradicts the fact that $\mathbb{R}^{m-m_1}$ is not compact.
Finally, $K$ is not positive.

Now let us prove that all eigenvalues of $L$ are constant. Without
loss of generality, we can assume that the manifold is simply
connected.  We will construct a totally geodesic (immersed)
submanifold $M_A$, which is a global analog of the submanifold
$M_A$  from Section~\ref{Vk}. At every point $x\in M^n$, consider
$V_{m_1+1},...,V_{m}\subset T_xM^n$, where $V_{m_1+i}$ is the
eigenspace of the eigenvalue $\phi_{m_1+i}$. Since the eigenvalues
$\phi_{m_1+i}$ are isolated, $V_{m_1+1},...,V_{m}$ are
distributions. By Theorem~\ref{nijenhuis}, they are integrable.
Denote by $M_{m_1+1}, M_{m_1+2},...,   M_{m}$ the corresponding
integral submanifolds.

 Since $M^n$ is simply connected,  by
\cite{hebda},
  it is diffeomorphic to
 the product
 $M_0\times M_{m_1+1} \times M_{m_1+2} \times ... \times M_{m}$.
 Clearly, the metric $g$ on $$M^n\simeq M_0\times M_{m_1+1} \times M_{m_1+2} \times ... \times M_{m}$$
 has the form
\begin{equation}\label{gM0} g_{|M_0} +\left|\prod_{\begin{array}{c}i=1\\ i\ne m_1+1\end{array}}^{m}
 (\phi_{m_1+1}-\phi_i)\right|g_{m_1+1}+...+\left|\prod_{\begin{array}{c}i=1\\ i\ne m \end{array}}^{m}
 (\phi_{m}-\phi_i)\right|g_{m},\end{equation}
where every $g_k$ is a  metric on $M_k$. Take a point
$$P=(p_0,p_{m_1+1},p_{m_1+2},...,p_m)\in M_0\times M_{m_1+1} \times
M_{m_1+2} \times ... \times M_{m}.$$ On every $M_{m_1+k}$,
$k=1,...,m-m_1$,  pick a geodesic $\gamma_{m_1+k}$ (in the metric
$g_{m_1+k}$) passing through $p_k$. Denote by $M_A$ the product $$
M_0\times \gamma_{m_1+1}\times ... \times \gamma_{m}. $$ $M_A$ is
an immersed totally geodesic manifold. Hence, it is complete in the metric $g_{|M_A}$ and in the metric  $\bar g_{|M_A}$. 
 Locally, in   a neighborhood of every point, $M_A$
 coincides with $M_A$ from Section~\ref{Vk} constructed for the warped decomposition
 (\ref{gM0}). The
restriction of the metric $g$ to $M_A$ is isometric to (\ref{g0})
and, therefore, has nonpositive constant sectional curvature $K$.
Then, by Corollary~\ref{bel}, the restriction of $\bar g$ to $M_A$
is affine equivalent to the restriction of $g$ to $M_A$. Then, by
Remark~\ref{eisen}, all $\phi_i$ are constant. Then,  $g$ is
affine equivalent to  $\bar g$. Theorem~\ref{4.2} is proved.

\section{Proof for Dimension 2} \label{dim2}
Let $(M^2,g)$ be a complete connected Riemannian manifold.    Our
goal is to proof Theorem~\ref{main} in dimension two. In view of
Corollary~\ref{bel}, we should show that  the existence of a
one-parametric family of nonaffine projective transformations of
$(M^2,g)$ implies that $g$ has constant   curvature.

 If the dimension of ${\cal B}(M^2,g)$
is 1, Theorem~\ref{main} is trivial. If the  dimension of ${\cal
B}(M^2,g)$ is two, Theorem~\ref{main}  follows from
Theorem~\ref{2dim}.

By \cite{dissertation,Kio}, if the manifold is closed, and if the
dimension of ${\cal B}(M^2,g)$ is greater than two, the curvature
of $g$ is constant. Then, it is sufficient to prove
Theorem~\ref{main} assuming that $M^2$ is not compact.

There exist  examples of complete metrics of nonconstant curvature
on noncompact surfaces   such that the dimension of ${\cal B}$ is
greater than two. Here are the    most important for us examples.
The manifold is $\mathbb{R}^2$ with the standard coordinates
$(x,y)$ and  $\gamma>0$ is a constant. The metrics are:

\noindent
{\bf Example 1:} $(x^2+y^2+\gamma)(dx^2+dy^2)$,

\noindent {\bf Example 2:} $(x^2+\frac{1}{4}
y^2+\gamma)(dx^2+dy^2)$.

 Example 1 was constructed in  \cite{konigs}.
Example 2 seems to  be new. We will explain in
Section~\ref{proofexamples} why the dimension of ${\cal B}$ is
four for Example 1 and three for Example 2.

If $M^2$ is  not compact, its   universal cover  is (homeomorphic
to) the 2-disk $D^2$. Hence,  it is sufficient to prove
Theorem~\ref{main} assuming that $M^2$   is  $D^2$. We will prove
the following stronger

\begin{Th}[Announced in \cite{obata,lich}] \hspace{-2mm}{\bf .} \label{d2}
 Suppose complete  Riemannian metrics $g$ and $g_1$  on the 2-disk $D^2$
 are projectively equivalent. Let them be not proportional (i.e. let
 $g\ne \textrm{const} \ g_1$
 for every $\textrm{const}\in\mathbb{R}$).
 If the curvature of $g$ is not
 constant, the space ${\cal B}(D^2,g)$ is  two-dimensional.
\end{Th}

 This theorem follows immediately from
 Theorems~\ref{complete}, \ref{last},  \ref{examples}, which we formulate below  and prove in
 Sections~\ref{proofcomplete},\ref{prooflast},\ref{proofexamples}, respectively.

It is known that there exists a  complex coordinate $z=x+iy$ on
the disk such that the metric $g$  has the form  $f(z)dzd\bar z$.
 (Here $\bar z$ denotes the complex conjugate to $z$). \weg{ It is known that the disk
with a complex structure is biholomorph to the complex plane
$\mathbb{C}$ or to a disk in the complex plane. The next Theorem
shows that,  under the assumptions of Theorem~\ref{d2},   the disk
with the complex coordinate $z$ is biholomorph  to $\mathbb{C}$.}

\begin{Th}\hspace{-2mm}{\bf .} \label{complete}
Suppose complete
 Riemannian metrics $g$ and $g_1$ are projectively equivalent
 on the disk  $D^2$. Let them be not proportional (i.e. let $g\ne
 \textrm{const} \ g_1$
 for every $\textrm{const}\in\mathbb{R}$). Then, the disk $D^2$
 with the complex structure $z$ such that
 the metric $g$  has the form $f(z)dzd\bar z$ is biholomorph  to the
 whole $\mathbb{C}$.
\end{Th}

\begin{Th}\hspace{-2mm}{\bf .}
\label{last} Consider a Riemannian   metric $g$ of the form
$f(z)dzd\bar z$ on $\mathbb{C}$. Let the dimension of ${\cal
B}(\mathbb{C}^2,g)$ be greater than two. Then, either $g$ is flat,
or there exists a diffeomorphism $F:\mathbb{C}\to \mathbb{R}^2$
taking the metric $g$ to   the metric from Example 1 or to the
metric from Example 2 for the appropriate $\gamma>0$.
\end{Th}

\begin{Th}\hspace{-2mm}{\bf .}\label{examples}
 If a  complete Riemannian  metric $\bar g$ on  $\mathbb{R}^2$ is  projectively equivalent to
 the metric from Example~1  or Example 2, then it is proportional to this metric.
\end{Th}

\subsection{Instruments of the proof} \label{instr}

\noindent Proofs of Theorems
\ref{complete},\ref{last},\ref{examples} use techniques  from the
theory of quadratically integrable geodesic flows on surfaces. The
relation between quadratically integrable geodesic flows and
projectively equivalent metrics is given  by the following
two-dimensional version of Theorem~\ref{integrability}:

\begin{Cor}[\cite{doklady,quantum,ERA}]\hspace{-2mm}{\bf .} \label{15}
Riemannian
metrics  $g$ and $\bar g$ on $M^2$
  are projectively equivalent,
 if and only if   the function
\begin{equation} \label{140} \label{I}
I:=  \left(\frac{det(g)}{det(\bar g)}\right)^{\frac{2}{3}}\bar g(\xi, \xi)
\end{equation}  is an integral of the geodesic flow of $g$.
\end{Cor}

The goal of this section is to formulate in a convenient form
important for us
results from the theory of quadratically integrable geodesic flows. Local form of most these results is classical and
can  be found in \cite{Darboux} or \cite{Birkhoff}; the global
generalizations  are  due to Kolokoltsov \cite{Kol,dissertation},
Kiyohara et al \cite{Igarashi,Kio,memo} and Bolsinov et al
\cite{BMF,BF}.

Let $g$ be a Riemannian metric on the 2-disk $D^2$.
 Consider the  complex coordinate $z$ on the disk such that the metric has the
form  $f(z)dzd\bar z$.  (Here and up to the end $\bar z$,
$\bar p$,
 $\bar a$  etc.
will denote the complex conjugation).

We will  work on the cotangent bundle to $D^2$. The tangent and
the cotangent bundles will always be identified by $g$. We denote
by $p$ the corresponding (complex) coordinate on the cotangent
bundle $T^*D^2$. Consider the real-valued quadratic in velocities
function $$ I(z,p)=a(z)p^2+b(z)p\bar p +\bar a(z) \bar p^2. $$
($a$ and $b$ are not assumed to be holomorphic; moreover, since
$F$ is real-valued,    $b$ must be   real-valued as well.) If we
make a coordinate change $z=z(z_{new})$, the coefficients $a$ and
$ b$ will be changed as well; the next lemma controls  how do they
change. In this form, it is due to Kolokoltsov; in a less
convenient form,  it was known at least to Birkhoff
\cite{Birkhoff}. Consider $$A:=-\frac{1}{a(z)} dz\otimes dz.$$

 For the energy integral

\begin{equation}\label{energy}
\frac{p\bar p}{f(z)},
\end{equation}
 $A$ is not defined.

\begin{Lemma}[\cite{Kol,dissertation}]\hspace{-2mm}{\bf .}
\label{kol}  If   $I$ is an integral of the geodesic flow of $g$,
and if it is not the energy integral (\ref{energy})  multiplied by
 a constant, then
 $A$ is a meromorphic (2,0)-form.
\end{Lemma}
\weg{{\bf Explanation: } After the  holomorphic change
$z=z(z_{new})$ of the coordinate, the momentum $p$ changes s
follows: $p_{new}=pz'$, where $z'$ denotes the derivative
$\frac{dz(z_{new})}{dz_{new}}$. Then, the integral $F$ in the new
coordinates
   is
\begin{eqnarray*}
F &=& a(z)p^2+b(z)p\bar p +\bar a(z) \bar p^2 \\
   &=& a(z(z_{new}))\frac{p_{new}^2}{(z')^2}+b(z(z_{new}))
\frac{p_{new}\bar p_{new}}{z'\bar z'} +\bar a(z(z_{new}))
\frac{\bar p_{new}^2}{(\bar z')^2}\\
 &=&
a_{new}(z_{new})p_{new}^2+b_{new}(z_{new})p_{new}\bar p_{new}
 +\bar a_{new}(z_{new})\bar p_{new}^2,
\end{eqnarray*}
where $a_{new}(z_{new})=\frac{a(z(z_{new}))}{(z')^2}.$ We see that
$-\frac{1}{a(z)}$ changes as a coefficient of a $(2,0)-$form. By
Corollary~\ref{ordered12}, if $I$ is not the energy integral
(\ref{energy})  multiplied by
 a constant, $a\ne 0$ almost everywhere so
that $A$ is defined at almost every point.

The fact that ${a(z)}$ is holomorphic has been known at least to
Birkhoff \cite{Birkhoff}; one  obtains  it immediately considering
the Poisson bracket of the energy integral (\ref{energy}) and $F$,
see, for example, \ref{kol}.

}

Locally, in a neighborhood of every point $P\in D^2$ which  is not a pole of $A$,
 by local holomorphic change of the variable
$w=w(z)$, we can always make the form $A$ to look $dw\otimes dw$.
Indeed, under this assumption $a(P)\ne 0$,  and    the equation
\begin{equation} \label{kol1}
 -\frac{1}{a(z)} dz\otimes dz = dw\otimes dw
\end{equation}
 has a solution $w(z)=\int \frac{dz}{\sqrt{-a(z)}}.$
 In this new  coordinate $w$,  the metric and the integral have the following
  very
 nice form  (this is a folklore known at least to Birkhoff; a  proof can be
 found in \cite{Kol}):

 \begin{Lemma}[\cite{Birkhoff,Kol}]\hspace{-2mm}{\bf
 .} \label{bir}  Let $I$ be a quadratic in velocities integral for
 the geodesic flow of the metric  $g$ on $D^2$.
Suppose the  form  $A$ from Lemma~\ref{kol} is equal to $dz\otimes
dz$.
 Then, in the coordinates $x:=\Re(z)$, $y:=\Im(z)$, the
metric and the integral $I$ have the following form:
 \begin{equation}\label{liouville1}
 ds_g^2=(X(x)-Y(y))(dx^2+dy^2),
 \end{equation}
$$
 I=\frac{Y(y)p_x^2+X(x)p_y^2}{X(x)-Y(y)},
$$
where $X$ and
$Y$ are functions of one variable.
\end{Lemma}

If the integral $I$ in Lemma~\ref{bir}  is constructed from
projectively equivalent metric $g_1$ by formula  (\ref{140}), the
metric $g_1$ reads:
 \begin{equation}\label{liouville2}
 ds_{g_1}^2=\left(\frac{1}{Y(y)}-\frac{1}{X(x)}\right)\left(\frac{dx^2}{X(x)}+\frac{dy^2}{Y(y)}\right).
 \end{equation}

\begin{Rem}\hspace{-2mm}{\bf .} \label{ind}
The form from Lemma~\ref{kol} constructed for the  linear
combination $\alpha I+\beta E$, where $E$ is the energy integral
 (\ref{energy}) and  $\alpha\ne 0$, is equal to $\frac{1}{\alpha}A$, where $A$ is
the form constructed for $I$. In particular, it has the same
structure of poles.
\end{Rem}

 A complete description for  complete metrics  whose geodesic
 flows
 admit
quadratic in velocities integrals which are  not  linear
combinations of the energy integral   and the square of an
integral linear in velocities,  was obtained in \cite{Igarashi}.
We reformulate
 their  results  as Theorem~\ref{14}.

Consider the triple $(D^2, g, A)$, where $D^2$ is a $2-$disk, $g$
is a Riemannian metric on it,  and $A$ is a meromorphic
(2,0)-form. Two such triples $(D_1^2, g_1, A_1)$ and $(D_2^2, g_2,
A_2)$ are said
 to be {\bf isomorphic}, if there exists a
diffeomorphism $H:D_1^2\to D_2^2$ that takes $g_1$ to
$g_2$ and  $A_1$ to $A_2$.

\begin{Th}[\cite{Igarashi}]\hspace{-2mm}{\bf .} \label{14}
Let  $g$ be a complete Riemannian metric on $D^2$. Suppose $I$ is
a quadratic
 in velocities
 integral  of the geodesic flow of $g$. Assume in addition that it
  is not a linear combination
of the square of an integral linear in velocities and the energy
integral. Then, the triple $(D^2, g, A)$, where $A$ is the
meromorphic (2,0)-form from Lemma \ref{kol},  is isomorphic to one of
the following model triples (for the appropriate parameters $r_1,
R_1, r_2, R_2$  and  $f$).

{\bf Model 1a:} Parameters: $r_1<R_1,r_2<R_2\in
\mathbb{R}\cup-\infty\cup+\infty$; $f$ is a positive function on
the disk $$ D^2:=\{z\in \mathbb{C}: \ \ r_1<\Re(z)<R_1; \  r_2<\Im
(z)<R_2\}. $$ The metric $g$ is given by $f(z)dzd\bar z$; the form
$A$ is equal to $dz\otimes dz$.

{\bf Model 2:}  Parameters: $R_1,R_2\in \mathbb{R}_+\cup+\infty$; $f$ is
a positive function on the disk $$ D^2:=\{z\in \mathbb{C}: \ \
r_2>|\Re(2\sqrt{z})|; R_2>\Im(2\sqrt{z})\}. $$ The metric $g$ is
given by $f(z)dzd\bar z$; the form $A$ is equal to
$\frac{1}{z}dz\otimes dz$.

{\bf Model 3:}  Parameters: $R_1\in \mathbb{R}\cup+\infty$;
 $f$ is a positive function on
the disk
$$
 D^2:=\{z\in \mathbb{C}: \ \
R_1>\Im\left(\arcsin\bigl({z^2-1}\bigr)\right)\}.
 $$
 The metric $g$ is
given by $f(z)dzd\bar z$; the form $A$ is equal to
$\frac{1}{z^2-1}dz\otimes dz$.

\end{Th}

The case when the integral is a linear combination  of the square
of an  integral linear in velocities and of the energy integral is
much more easy than the previous case, and, may be because of it,
was not considered in \cite{Igarashi}. We will need it.

\begin{Th}\hspace{-2mm}{\bf .} \label{linear}
Let $g$ be  a complete metric  on $D^2$. Suppose
  the function $F$ is  a
linear combination of the square of an   integral linear in
velocities and the energy integral and is not proportional to the
energy integral. Then, the triple $(D^2, g, A)$, where $A$ is the
meromorphic (2,0)-form from Lemma~\ref{kol}, is isomorphic to one
of the following model triples (for the appropriate parameters $r,
R$ and $f$).

{\bf Model 1b:} Parameters: $r<R\in \mathbb{R}\cup+\infty\cup -\infty$;
$f$ is a positive function depending only on the variable
$x:=\Re(z)$  on the disk $$ D^2:=\{z\in \mathbb{C}: \ \
r<\Re(z)<R; \}. $$ The metric $g$ is given by $f(z)dzd\bar z$; the
form $A$ is equal to $dz\otimes dz$.

{\bf Model 4:}  Parameters: $R\in \mathbb{R}\cup+\infty$;
 $f$ is
a positive function depending only on the absolute value of $z$
on the disk $$ D^2:=\{z\in \mathbb{C}: \ \ |z|<e^R \}. $$ The
metric $g$ is given by $f(z)dzd\bar z$; the form $A$ is equal to
$\frac{1}{z^2}dz\otimes dz$.

\end{Th}
{\bf Proof:} Suppose the linear in momenta function
$I=a(z)p_x+b(z)p_y$ is an integral for the geodesic flow of $g$.
Then,  $v=(a,b)$ is a Killing vector field for $g$, and its flow
preserves the metric.
 Suppose first there exists no point where $v$ vanishes. Consider
 the vector field $w$ such that
 \begin{enumerate}
\item $g(v,v)=g(w,w);$ \ $g(v,w)=0$. \label{w}
\item The pair of the vectors $w,v$ is positively oriented.
 \end{enumerate}
These two conditions define the vector field $w$ uniquely.
Since $w$ is defined using $g$ and $v$ only, and since $g$ is preserved by the flow of $v$,
 the vector fields $w$ and $v$ commute. Then,   they  define  a coordinate system
$(x,y)$ on  the disk. We treat $w$ as the first vector field, so
that
the metric is independent of the variable $y$.  By the condition
\ref{w}, the metrics  has the form $f(z)(dx^2+dy^2)$.  Then,  $f$
depends on $x$ only. Since  the metric is complete,
the  lines $\{(x,y)\in D^2:\ x=\textrm{Const}\}$ are infinite in
both directions.

Then, there exist $r< R \in
\mathbb{R}\cup +\infty \cup -\infty$ such that the disk is
actually the band
$$
\{(x,y)\in D^2:  \ \ r< x
<R.\}
$$
By construction,  $v=(0,1)$. Then,
  the square of the integral $I$ is equal to $p_y^2$,  and
 the  form $A$ is $dz\otimes dz$.
  Theorem~\ref{linear} is proved
under the assumption that  $v$ vanishes at no points.

Now suppose there exist points where $v$ vanishes. These points
are  stable points of the flow of $v$; since the flow acts by
isometries, they are isolated.  Let us show that there exists
precisely one such point. Indeed, if $v$ vanishes at the points
$z_0\ne z_1$, since our manifold is homeomorphic to the disk and
is complete, there exist two geodesics $\gamma_0$, $\gamma_1$
 such
that
\begin{itemize}
\item The geodesic $\gamma_0$ contains $z_0$;   the geodesic $\gamma_1$ contains
$z_1$.

\item The geodesics intersect transversally at a point $z$ where $v\ne
0$.
\end{itemize}
Since the geodesics pass through  the point where $v$ vanishes,
the integral $I$ is equal to zero on  every velocity vector of the
geodesics, and, therefore, the geodesics are orthogonal to $v$ at
every point. This gives us a contradiction at the point $z$. Thus,
$v$ vanishes precisely at one point. We will denote this point by
$z_0$.

Since
  the flow of the field
$v$ acts by isometries and preserves the point $z_0$, it commutes
 with the exponent mapping $\exp_{z_0}:T_{z_0}D^2\to
D^2$. Since $v$ is a Jacobi vector field for geodesics passing
through $z_0$, there is no point conjugate to $z_0$ and,
therefore, the exponential mapping is a bijection. The flow of $v$
acts  on $T_{z_0}D^2$  by  orthogonal linear transformations,
i.e., by  rotations.
 Thus, in the standard polar  coordinates on $T_{z_0}D^2$,
  the pull-back of the metric is $d\rho^2+f_1(\rho)d\phi^2$, and
  $v$ is proportional to $\frac{\partial }{\partial \phi}$.
After the appropriate change of the variables, the metric is
  $f(|z|)dz\bar dz$,
 and the vector field $v$  is (proportional to)
 $(-y,x)$. Then,  the square of the integral $I$  is
 $y^2p_x^2-2xyp_xp_y+x^2p^2_y$, and the form $A$ is (after the scaling)
 $\frac{1}{z^2}dz\otimes dz$.
  Theorem~\ref{linear} is proved.

\begin{Rem}\hspace{-2mm}{\bf .}
Clearly, not arbitrary parameters $(R_i,r_i,f)$ can come from
metrics with quadratically integrable geodesic flows: for example,
the metric $g$ must have the form (\ref{liouville1}) in the
coordinates where $A$ is $dz\otimes dz$.
\end{Rem}

\subsection{Proof of Theorem~\ref{complete}} \label{proofcomplete}
 Let $g$, $g_1$ be complete projectively equivalent
nonproportional  Riemannian metrics on $D^2$. Our goal is to prove
Theorem~\ref{complete}. That is, we have to prove that the complex
structure of the disk corresponding to the metric $g$ is as of
$\mathbb{C}$.

 Consider the integral (\ref{140}) for the geodesic
flow of $g$. It is quadratic in velocities and linear independent
of the energy integral (\ref{energy}). Denote by $A$ the form from
Lemma~\ref{kol}.  By Theorem~\ref{14}, the triple $(D^2,g,A)$ is
isomorphic to one of the model triples from
Theorems~\ref{14},\ref{linear}. We will consider all four cases.
Suppose $(D^2,g,A)$ is as in model 1. We need to prove that
$R_1,R_2$ are $+\infty$ and $r_1,r_2$ are $-\infty$. We will prove
that $R_1$ is $+\infty$. By Lemma~\ref{bir}, in the coordinates
$(x:=\Re(z), y:=\Im(z))$,  the metrics $g$ and $g_1$    have  the
form (\ref{liouville1},\ref{liouville2}), respectively. Take $t_1$
such that  $r_1<t_1<R_1$ and  consider the line $\{x+iy\in D^2: \
y=0; \ x>t_1 \}$. Since the metrics are complete, the length of
the line must be infinite in both metrics. The length in $g_1$ is
 \begin{equation}\label{ii}
 \int_{t_1}^{R_1}\frac{\sqrt{X(t)-Y(0)}}{X(t)\sqrt{Y(0)}}dt.
 \end{equation}
Since $X(t)>Y(0)$ by Corollary~\ref{ordered12}, and since
$Y(0)>0$,  we have $$
\frac{\sqrt{X(t)-Y(0)}}{X(t)\sqrt{Y(0)}}<\frac{1}{\sqrt{Y(0){X(t)}}}<
\frac{1}{{Y(0)}}. $$ Thus, if (\ref{ii}) is infinite, $R_1$ is
$+\infty$. Similarly, one can prove that $r_1,r_2=-\infty$ and
$R_2=+\infty$. Thus, $D^2$ is the whole $\mathbb{C}$.

Now suppose $(D^2,g,A)$ is as in model 2. We need to prove that
$R_1=R_2=+\infty$. Take $0<t_1<R_1$ and consider the part of the
disk $$\{z\in D^2\subset \mathbb{C}: \ \Re(2\sqrt{z})>t_1\}.$$ In
this part of the disk,  the equation (\ref{kol1}) has no critical
points and can be explicitly solved; the  solution is
$w=2\sqrt{z}$. After the substitution $z=\frac{w^2}{4}$, the part
of the disk becomes to be  the rectangle $$\{w\in \mathbb{C}:\
t_1<\Re(w)<R_1\}. $$ By Lemma~\ref{bir}, the metrics $g$, $g_1$ in
the coordinates $x_{\mathrm{new}}:=\Re(w),
y_{\textrm{new}}:=\Im(w)$ have the form
(\ref{liouville1},\ref{liouville2}), respectively. Arguing as for
model 1,  since  the metric $ g_1$ is complete, (\ref{ii}) must be
infinite. Hence,  $R_1$ is $+\infty$. Similarly, one can prove
that $R_2=+\infty$. Thus, $D^2$ is the whole $\mathbb{C}$.

 Suppose $(D^2,g,A)$ is as in model 3. We need to prove that
$R_1=+\infty$. Take $t_1\in ]0,R_1[$ and consider the part of the
disk
$$
\{z \in D^2\subset\mathbb{C}:
t_1<\Im(\arcsin({z^2-1}))<R_1\ \}.
$$
 This part is homeomorphic
to the cylinder. The equation (\ref{kol1}) can be explicitly
solved on the universal cover of the cylinder; the solution is
$w=\arcsin({z^2-1}))$.    By Lemma~\ref{bir}, the metrics $g$,
$g_1$ in the coordinates $x_{\textrm{new}}:=\Re(w),
y_{\textrm{new}}:=\Im(w)$ have the form
(\ref{liouville1},\ref{liouville2}), respectively. Arguing as
above, if  $g_1$ is complete, $R_1=+\infty$.   Thus, $D^2$ is the
whole $\mathbb{C}$.

The last case is when  $(D^2,g,A)$ is as in model 4. We need to
prove that $R=+\infty$. Take $t_1\in ]0,R[$ and consider the part
of the disk $\{z \in D^2\subset\mathbb{C}: t_1<\Re(\log(z))<R\
\}.$ This part is homeomorphic to the cylinder. The equation
(\ref{kol1}) can be explicitly solved on the universal cover of
the cylinder; the solution is $w=\log(z)$.    By Lemma~\ref{bir},
the metrics $g$, $g_1$ in the coordinates
$x_{\textrm{new}}:=\Re(w), y_{\textrm{new}}:=\Im(w)$
 have the form
(\ref{liouville1},\ref{liouville2}), respectively; moreover, the
function $Y$ is constant.  If the metric $g_1$ is complete, the
length of the line $$ \{w\in \mathbb{C}: \ y_{\textrm{new}}=0, \
x_{\textrm{new}}\in]t_1,R[\} $$ must be infinite. Thus,
$\int_{t_1}^{R}\frac{\sqrt{X(t)-Y}}{X(t)\sqrt{Y}}dt$ is infinite,
so that $R$ is infinite.  Thus, $D^2$ is the whole $\mathbb{C}$.
Theorem~\ref{complete} is proved.

\subsection{Proof of Theorem~\ref{last}} \label{prooflast}
We will need the following corollary of
Theorems~\ref{14},\ref{linear}.
\begin{Cor}\hspace{-2mm}{\bf .} \label{1234}
Consider  a  complete Riemannian metric $g$ on $\mathbb{C}$ given
by $ds_g^2=f(z)dzd\bar z$.  Suppose the  geodesic flow has an
integral $I$ that is quadratic in velocity and functionally
independent of the energy integral (\ref{energy}).
  Denote by $A$ the
 (2,0)-form from Lemma~\ref{kol}. Then, the
  following statements hold:
\begin{enumerate}
 \item The form $A$ is  $\frac{1}{\alpha z^2+\beta z+ \gamma}dz\otimes
 dz$, where $\alpha$,  $\beta $ and  $\gamma$ are complex
 constants. At least one of these constants is not zero.

  \item If  $\alpha=0$  and
  $\beta\ne 0$, the metric is
  preserved by the symmetry with respect to the straight line
  $$
\{z\in \mathbb{C}:\ \frac{z}{\beta }  + \frac{\gamma}{\beta^2} \in
\mathbb{R}\}.
  $$

 \item If the polynomial  ${\alpha z^2+\beta z+ \gamma}$ has two
 simple roots, the metric  $g$  is
 preserved by  the symmetry with respect to the straight line
 connecting the roots.

 \item If the polynomial  ${\alpha z^2+\beta z+ \gamma}$ has one double  root,
  the integral (\ref{I}) is a linear combination of the  energy
  integral (\ref{energy})  and the square of an  integral linear in velocities,
  which vanishes at the tangent plane to the root. The rotations
  around the root preserve the metric.

\end{enumerate}
\end{Cor}
\noindent{\bf Proof:}
By Theorem~\ref{14},
there exists complex coordinate  $w$
 on $\mathbb{C}$ such that $A$ is
either $\frac{1}{C}dw\otimes dw$ or  $\frac{1}{w}dw\otimes dw$ or
$\frac{1}{w^2-1}dw\otimes dw$ or $\frac{1}{w^2}dw\otimes dw$.
Since every bijective holomorphic mapping $:\mathbb{C}\to
\mathbb{C}$ is linear, the change $w=w(z)$  of the coordinate  is
linear. After a linear change of the coordinate, the forms $A$
listed in Theorem~\ref{14},\ref{linear},   have  the form
$\frac{1}{\alpha z^2+\beta z+ \gamma}dz\otimes
 dz$, where $\alpha$,  $\beta $ and  $\gamma$ are complex
 constants such that at least one of them is not zero.
 The first statement of Corollary~\ref{1234} is proved.

 In order to prove the second statement, it is sufficient to show
 that the symmetry $z\mapsto \bar z$ is an isometry of the metric
 whose triple is $(\mathbb{C}^2,f(z)dzd\bar z,\frac{1}{z}dz\otimes dz)$.
 Indeed, after  the
 following change of coordinate
 $$
 w=\frac{z}{\beta}+\frac{\gamma}{\beta^2},
 $$
 the form $A$ is $\frac{1}{w}dw\otimes dw$ and the line
  $
\{z\in \mathbb{C}: \ \frac{z}{\beta }  + \frac{\gamma}{\beta^2} \in
\mathbb{R}\}
  $
  becomes the line
 $\{w\in\mathbb{C}:\ \Im(w)=0\}$.

\begin{figure}[p!]
  \hspace{-2cm}
{{\psfig{figure=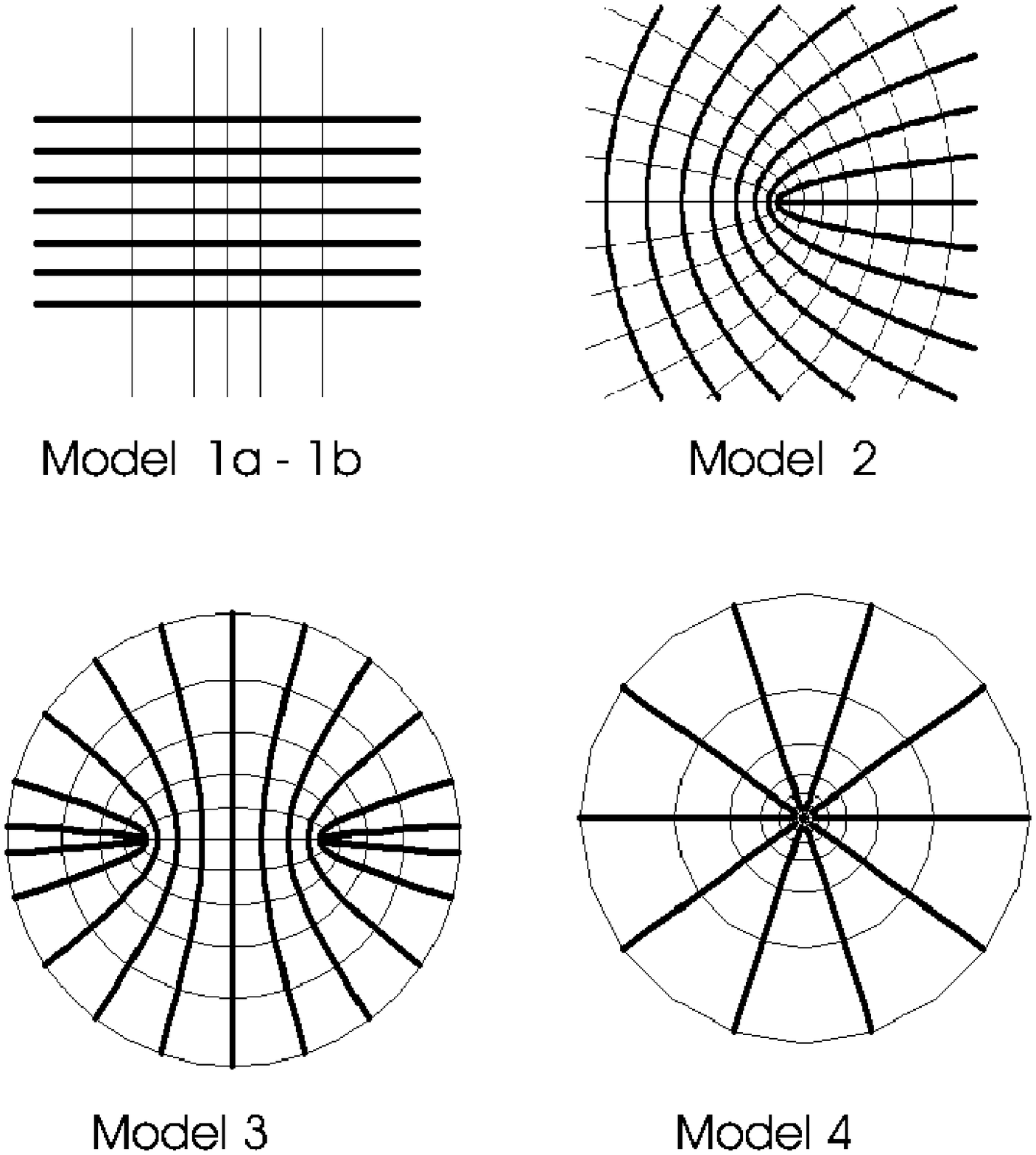}}}
  \caption{Coordinates where $A=dz\otimes dz$ for models 1--4}\label{picture}
\end{figure}

 Clearly, the form $A=\frac{1}{z}dz\otimes dz$ has precisely one pole. In
 the neighborhood of every  other point, the equation
 (\ref{kol1}) can be explicitly  solved: the solution is
 $w=2\sqrt{z}$. It can not be solved globally: the global solution is
 defined on the double branch cover. But still the coordinate
 lines of the local coordinate system
 $(\Re(2\sqrt{z}),\Im(2\sqrt{z}))$ do  not depend on the choice of
 the branch of the square root and can be defined globally. These lines look as on the
 Figure~\ref{picture}, model 2.
 We see, that if a
 coordinate line of $\Re(w)$ intersects with a coordinate
 line of $\Im(w)$ at a point $z$, then the same coordinate lines
 intersect at the point $\bar z$ as well. By Lemma~\ref{bir}, the
 metric $g$ has the form
 $$
\frac{1}{\sqrt{z\bar z}}\left(X(\Re(w))-Y(\Im(w))\right)dzd\bar z.
$$
 Since the function $X$ is constant along the coordinate lines of
 $\Im(w)$ and the function $Y$ is constant along the coordinate lines of
 $\Re(w)$, $X$ at  $z$ is equal to $X$ at  $\bar z$ and
 $Y$ at  $z$ is equal to $Y$ at  $\bar z$. Thus, the
 symmetry $z\mapsto \bar z$ is an isometry.
 Statement 2 is proved.

 The proof of  statement 3 is similar: the solution of
 Equation (\ref{kol1}) for the model triple 3 is
 $w=\arcsin({z^2-1})$, and the coordinate lines such that
  the
 metric has the form (\ref{liouville1}) are as on Figure~\ref{picture}.

 Fourth statement follows from
 Theorems~\ref{14}, \ref{linear}. Corollary~\ref{1234} is proved.

\thispagestyle{empty}

\weg{
\hspace{4cm}\fbox{\includegraphics{case_small1.eps}}
}

{\noindent \bf Proof of Theorem~\ref{last}:} Consider two
quadratic in velocities integrals  $I_1,I_2$ such that $I_1,I_2$
and the energy integral (\ref{energy}) are linear independent.
Denote by $A_1$ (respectively, by $A_2$)  the (2,0)-form from
Lemma~\ref{kol} constructed for the integral $I_1$ (respectively,
$I_2$).

Every   triple $(\mathbb{C},g,A_i)$ must be isomorphic to one of
the four model triples. Below we will consider all ten possible
cases.

\noindent{\bf Case  (1,1):} Suppose    both triples
 are isomorphic
to   model  1. Then, after the appropriate change of coordinates,
the form $A_1$ is $dz\otimes dz$ and the form $A_2$ is
$C^2dz\otimes dz$, where $C=\alpha+\beta i$ is a complex constant.
Then, the metric and the integrals have the model form  from
Lemma~\ref{bir} in the coordinate systems $(x:=\Re(z),y:=\Im(z))$
and $\bigl(x_{\textrm{new}}:=\alpha x- \beta y=\Re (Cz),
y_{\textrm{new}}:=\alpha y+\beta x=\Im (Cz)\bigr).$ Then, for
certain functions $X_1,X_2,Y_1,Y_2$ of one variable,
 $$
ds^2_g=(X_1(x)-Y_1(y))(dx^2+dy^2)=(X_2(\alpha x- \beta y)-Y_2(\alpha
y+\beta x))(dx^2+dy^2). $$

If $\alpha$ or $\beta$ is zero, the integrals $I_1,I_2,E$ are
linear dependent, which contradicts the assumptions. Assume
$\alpha\ne 0 \ne \beta$.  Since $ \frac{\partial^2
(X_1(x)-Y_1(y))}{\partial x
\partial y}=0, $ we have $$ \frac{\partial^2
 }{\partial x \partial y}(X_2(\alpha x-
\beta y)-Y_2(\alpha y+\beta x))= -\alpha\beta(X''_2(\alpha x-
\beta y) + Y''_2(\alpha y+\beta x))=0. $$ Thus, $X''_2(\alpha x-
\beta y) = - Y''_2(\alpha y+\beta x))=\textrm{Const}.$
 Finally, $X_2=-Y_2=\textrm{Const}_1$,  or
$X_2$ and $-Y_2$ are   quadratic  polynomials   with the same
coefficient near the quadratic term,  or  $X_2$ and $-Y_2$  are
 linear polynomials. In the  first case, the
metric is flat.  In the second   case,  after the appropriate
change of variables, the metric  is as in Example 1. In the third
case, the metric is not positive definite. Theorem~\ref{last} is
proved under the assumption that both triples are as in  model 1.

\noindent{ \bf Case (1,4):} Suppose
 $A_1=\frac{1}{C}dz\otimes dz, \ \ A_2=\frac{1}{z^2}dz\otimes dz.$
 Then,
  the metric is
invariant with respect to the rotations $z\mapsto e^{i\phi}z$.
Then, the push-forward of the integral is an integral.  By
Remark~\ref{ind},  for $\phi\not\in \frac{\pi}{2}\mathbb{Z}$, the
push-forward of the integral $I_1$ is linear independent of the
integral $I_1$ and of the energy integral. Its form from
Lemma~\ref{kol} has no poles. Thus, we reduced
 case (1,4) to  case (1,1).

\noindent{\bf Case  (1,2):} Suppose
 $A_1=dz\otimes dz, \ \ A_2=\frac{1}{\beta z+\gamma}dz\otimes dz.$
      Then,
 by Lemma~\ref{bir},
 the metric $g$ is
$(X(x)-Y(y))(dx^2+dy^2)$.
    If $\beta $ is not a real multiple  of $1,
i$ or $1\pm i$, by Remark~\ref{ind}, the isometry from
Corollary~\ref{1234} sends the integral $I_1$ to the integral that
is not a linear combination of the integral $I_2$ and the energy
integral. Clearly, the form from Lemma~\ref{kol}  for this
integral has no pole, so we reduced this case to case~(1,1).

Suppose  $\beta$ is a real multiple of  $1+ i$. Without loss of
generality, we can assume  $\beta=1+ i$. Consider the integral
$tI_1+I$. Its form from Lemma~\ref{kol} is $$ \frac{1}{(1+i) z+
t+\gamma}. $$ Then, by Corollary~\ref{1234}, for every $T\in
\mathbb{R}$, the symmetry with respect to the line $$ \left\{z\in
\mathbb{C}: \frac{z}{1+i} + \frac{\gamma+T}{2i} \in
\mathbb{R}\right\} $$ is an isometry of $g$.  This one-parameter
family of symmetries gives us a Killing vector field $v=(1,1)$.
Thus, for every $z=x+iy$ and for every real constant $c$,
$X(x+c)-Y(y+c)=X(x)-Y(y)$. Thus, $X$ and  $Y$ are linear functions
or constants. If they are constants, the metric is flat. If they
are linear functions, the metric is not everywhere  positive
definite. Thus, $\beta$ is not a real multiple of  $1+ i$.

The proof that   $\beta$ is not  a real multiple of  $1- i$ is
similar.

Now suppose $\beta = 1$ or $\beta = i$.  Without loss of
generality, we can assume  $\gamma=0$.

Then, by Lemma~\ref{bir},    the metric has the Liouville  form (\ref{liouville1}) in  coordinate systems $\bigl(\Re(z), \Im(z)\bigr)$ and $\bigl(\Re(2\sqrt{z}), \Im(2\sqrt{z})\bigr)$, which imply
$$
ds^2=\left(X_1\bigl(\Re(z)\bigr)-Y_1\bigl(\Im(z)\bigr)\right)dzd\bar z = \left(
X_2\bigl(\Re(w)\bigr)-Y_2\bigl(\Im(w)\bigr)\right)dw d\bar w,
$$
 where
$w=2\sqrt{z}$. We see that   $dz$ and $dw$ are connected by the
relation
\begin{equation}\label{-4}
dw=\frac{1}{\sqrt{z}}dz,\end{equation} and
$(X_2-Y_2)=|z|(X_1-Y_1).$\weg{ By Lemma~\ref{bir}, the metric  has the
form ~(\ref{liouville1}) in the coordinates $(\Re(w),\Im(w))$.}
Then,
 the imaginary part of
$$ \frac{\partial^2 (X_1(x)-Y_1(y)) |z|}{\partial w^2} $$ must be
zero. Substituting  (\ref{-4}), we obtain
\begin{equation}
X_1''y+Y_1''y+3Y_1'=0.
\end{equation}
This equation can be solved. The solutions are $$
X_1(x)=C_1x^2+C_2x+C_3, \ \
Y_1(y)=-\frac{1}{4}C_1y^2+C_4\frac{1}{{y^2}}+C_5. $$  Since $X_1$
 and $Y_1$ must be continuous, $C_4=0$.  If $C_1=0$, the metric is flat or not positive definite.
If $C_1\ne 0$, after the appropriate linear change of coordinates,
the metric is as in Example 2. Theorem~\ref{last} is proved under
the assumption that one integral is as in model 1, and the other
is as in model 2.

 \noindent{ \bf Case (1,3):}
Suppose  $A_1=\frac{1}{C}dz\otimes dz$ and  $A_2=\frac{1}{z^2-1}
dz\otimes dz$. Then, the form from Lemma~\ref{kol} corresponding
to the integral $\alpha I_1 + I_2$ is $\frac{1}{z^2-1 +\alpha
C}dz\otimes dz$. If $C$ is real, the form for   $\frac{1}{C} I_1 +
I_2$  is as in model 4, i.e. we can  reduce case (1,3) to case
(1,4).

Suppose $C$ is not real.  Then, by Corollary~\ref{1234}, for every
real $\alpha$, the metric $g$ is preserved by the symmetry with
respect to the line connecting $\pm\sqrt{1-\alpha C}$. Then, the
rotations $z\mapsto e^{i\phi}z$ are isometries. We reduced again case
(1,3) to case (1,4).

\noindent{ \bf Case (2,2):} Suppose $A_1=\frac{1}{z}dz\otimes dz$
and $A_2=\frac{1}{\beta z+\gamma }dz\otimes dz$. If $\beta$ is
real,  the form from Lemma~\ref{kol} for  $-\beta I_1+I_2$ has no
pole, so that we reduced case (2,2) to  case (1,2). If $\beta$ is
not real, the form from Lemma~\ref{kol} for integral $t I_1 + I_2$
is $\frac{1}{(\beta -t)z+\gamma }dz\otimes dz$. We see that the
line of the symmetry from Corollary~\ref{1234} smoothly depend on
$t$, and is not constant. Then, the symmetries from
Corollary~\ref{1234} generate an one-parametric family of
isometries of $g$. By Noether's Theorem, a family of isometries
generates an integral linear in velocities. We consider the square
of this linear integral. By Theorem~\ref{linear}, either there
exists a point such that the integral
 vanishes at the tangent space to the point, or
there is no such point. In the second case, we reduced case (2,2)
to  case (1,2). In the first case, if the point where the vector
field vanishes does not lie on the line $\Re(z)=0$, the symmetry
of the integral  with respect to the line is also an integral, so
that  we constructed two linear independent Killing vector fields.
Then, the metric is flat.

The only remaining possibility  is when the point where the vector
field vanishes lies on the line $\Re(z)=0$. If the point does not
coincide with the point $0$,  the symmetry with respect to the
point gives us another integral with the form from Lemma~\ref{kol}
equal to $\frac{1}{C-z} dz\otimes dz$. Then,  a linear combination
of this other  integral and the integral $I_1$ gives us an integral such
that the  form from Lemma~\ref{kol} has no pole. Thus, we reduced
case (2,2) to case (1,2).

Now suppose the point where the Killing vector field vanishes
coincides with $0$. Then, by Theorem~\ref{linear}, the metric has
the form (we assume $w=2\sqrt{z}$)  $$ f(\sqrt{z\bar z})dzd\bar
z=\frac{1}{4}w\bar w f\left(\frac{w\bar w}{4}\right ) dwd\bar w
.$$ Since by Lemma~\ref{bir}  the metric has the
form~(\ref{liouville1}) in the coordinates $\Re(w),$ $\Im(w)$, we
obtain $\frac{1}{4}w\bar w f(\frac{w\bar w}{4})=C_1w\bar w +C_2$,
which imply that the metric is
 $ {}\left(4C_1    +\frac{C_2}{\sqrt{z\bar z}}\right)dz d\bar z$. We see
 that the metric is  flat or degenerate.
Theorem~\ref{last} is proved under assumptions of case~(2,2).

\noindent{ \bf Case (2,4): }
 Suppose $A_1=\frac{1}{\beta z+ \gamma}dz\otimes dz$
and $A_2=\frac{1}{z^2}dz\otimes dz$. Then, the metric is invariant
with respect to the
 rotations $z\mapsto e^{i\phi}z$. Arguing as in case (1,4),
  we can construct   one more integral
 such that it is linear independent
of $I_1$ and $E$, and such that the  form form Lemma~\ref{kol} is
as in model 2. Thus, we reduced case~(2,3) to case (2,2).

\noindent{ \bf Case (2,3): }  Suppose $A_1=\frac{1}{z}dz\otimes dz$
and $A_2=\frac{1}{\alpha z^2+\beta z +\gamma  }dz\otimes dz$. The form
from Lemma~\ref{kol} for the linear combination $I_2+tI_1$ of the
integrals  is $\frac{1}{\alpha z^2+(\beta+t) z +\gamma}dz\otimes dz$.
Consider the symmetries from Corollary~\ref{1234} for the linear
combination $I_2+tI_1$ of the integrals. We know that they are
symmetries with respect to the line connecting the roots of the
polynomial $\alpha z^2+(\beta+t) z +\gamma$. Analyzing these
symmetries, we see that they do not generate  one more integral
corresponding to the model 2  such that it is linear independent
of  $I_1$ and $E$, if and only if $\alpha$, $\beta$ and  $\gamma$ are
 real. In this case, a linear combination has a double
root.
 Thus, we reduced case~(2,3) to  cases (2,2),  (2,4).

\noindent {\bf Case (4,4):} In this case,  we have a  two-parameter
group of isometries, which is possible only if the metric is flat.

\noindent{\bf Case (3,4):} Suppose $A_1=\frac{1}{z^2} dz\otimes
dz$ and $A_2=\frac{1}{\alpha z^2+\beta z+\gamma} dz\otimes dz$.
Then,  the metric is invariant with respect to the rotations
around $0$. Consider the line connecting the roots of $\alpha
z^2+\beta z+\gamma$. If $0$ does not lie on the line, the point
symmetric to $0$ with respect to this line does not coincide with
$0$, and we reduce case (3,4) to case (4,4).

Now suppose $0$  lies  on the line connecting the roots
 of $\alpha z^2+\beta z+\gamma $.   Then,
a linear combination of  the square of the linear integral  and of
integral $I$ is as in model 2; so the we reduce case (3,3) to case
(2,3).

\noindent{ \bf Case (3,3):} Suppose $A_1=\frac{1}{z^2-1} dz\otimes dz$
and $A_2=\frac{1}{\alpha z^2+\beta z+\gamma } dz\otimes dz$.
 Consider the symmetries from
Corollary~\ref{1234} for the linear combination $I_2+tI_1$ of the
integrals. It is easy to see that if not all these symmetries
coincide, they generate  at least  one-parametric family of
isometries of $g$. By Noether's Theorem, this family generates an
integral linear in velocities of the geodesic flow.  The square of
the integral is as in model 1b (so we reduce case (3,3) to case
(1,3)), or  as in model 4 (so we reduce case (3,3) to case (3,4)).
If all the symmetries coincide, a linear combination of the
integrals has a double root, so we reduce case (3,3) to case
(3,4).

Finally, in all 10 cases we considered, the metric is  either as in
Examples 1,2, or  flat. Theorem~\ref{last} is proved.

\subsection{Proof of Theorem~\ref{examples}}\label{proofexamples}

Our goal is to show that if a complete Riemannian metric is
projectively equivalent to a  metric from Examples~1,2, then it is
proportional to the metric. Let us do it for Example 1. First of
all, the space ${\cal B}(\mathbb{R}^2,g)$ for the metric
$(x^2+y^2+\gamma)(dx^2+dy^2)$ on $\mathbb{R}^2$ has dimension
four. Indeed,  we can present four linear independent integrals.
They are (in  the standard coordinates $(x,y,p_x,p_y)$ on
$T^*\mathbb{R}^2$):

\begin{eqnarray*}
 H&:=&{\frac {{{\it p_x}}^{2}+{{\it p_y}}^{2}}{{x}^{2}+{y}^{2}+\gamma}},\\
 F_1&:=&{\frac {{{\it p_x}}^{2}{y}^{2}- \left( {x}^{2}+\gamma \right) {{\it p_y}
}^{2}}{{x}^{2}+{y}^{2}+\gamma}}\\ F_2&:=&(x{\it p_y}-y{\it
p_x})^2\\ F_3&:=&  xyH-p_xp_y
\end{eqnarray*}

The integrals are clearly linear independent (for example because
  the form from Lemma~\ref{kol} for  a nontrivial linear combination of
the   integrals $H$  and $F_1$ is $C_1\, dz\otimes d z$ and is
never equal to  the
   form
from Lemma~\ref{kol} for  a nontrivial  linear combination of
  the integrals $F_2$ and $F_3$  which is
 $\frac{1}{C_2z^2+ C_3i} dz\otimes d  z$).

Since the curvature of  $g$ is not constant,  the dimension of
${\cal B}(\mathbb{R}^2,g)$ can not be greater than four by
\cite{konigs}.

Hence, every  integral quadratic in velocities is a linear
combination of the integrals above. Clearly, the poles of the form
from Lemma~\ref{kol}  for the linear combination of the integrals
are  symmetric with respect to the point $0$. Then,   modulo
rotation of the coordinate system and scaling,  we can assume that
the form is either $dz\otimes d z$ or $\frac{1}{z^2}dz\otimes d z$
or $\frac{1}{z^2-1}dz\otimes d z$. In the first case, the
projectively equivalent metric constructed by Corollary~\ref{15}
from the integral is $$
\pm\left(\frac{1}{C-y^2}-\frac1{x^2+\gamma+C}\right)\left(\frac{dx^2}{x^2+\gamma+C}+\frac{dy^2}{C-y^2}
\right),$$
 and we see that it is not always positive defined. In the
second case, the projectively equivalent metric  constructed from
the integral has the form $$\pm
\left(\frac{1}{C}-\frac{1}{(r^2+\gamma)r^2+C}\right)\left(\frac{dr^2}{r^2\left
((r^2+\gamma)r^2+C \right) }+\frac{dy^2}{C} \right)$$  in the
polar coordinates.  We see that the metric is not positive
definite or not  complete.

In the third  case, the projectively equivalent metric constructed
from the integral has the form $$
\pm\left(\frac{1}{X}-\frac{1}{Y}\right)\left(\frac{dx^2}{X}+\frac{dy^2}{Y}
\right)$$ in the coordinates $x_{\textrm{new}}+i
y_{\textrm{new}}:=\arcsin(x+iy)$,  where  the functions $X,Y$ are
given by
\begin{eqnarray*}
2X(x_{\textrm{new}})& = &  \frac{1}{2}
\left(\cos(2x_{\textrm{new}})\right)^2-\gamma\,
\cos(2x_{\textrm{new}})+C
\\ 2Y(y_{\textrm{new}})& = & \frac{1}{2}
\left(\cosh(2y_{\textrm{new}})\right)^2+\gamma\,
\cosh(2y_{\textrm{new}})+C.
\end{eqnarray*}
 We see that the metric is not positive definite or not complete.

Thus, every complete Riemannian  metric which is projectively
equivalent to the metric from Example 1 is proportional to it.

Now let us show that the metric from Example 2 admits no
nonproportional projectively equivalent Riemannian  metric.

The   space ${\cal B}(\mathbb{R}^2,g)$, where  $g$  is the metric
from Example 2,  has dimension three. Indeed,  we can present
three linear independent integrals. They are (in  the standard
coordinates $(x,y,p_x,p_y)$ on $T^*\mathbb{R}^2$):

\begin{eqnarray*}
 H&:=&{\frac {{{\it p_x}}^{2}+{{\it p_y}}^{2}}{{x}^{2}+\frac{1}{4}{y}^{2}+\gamma}},\\
 F_1&:=&{\frac {{{\it p_x}}^{2}\frac{{y}^{2}}{4}- \left( {x}^{2}+\gamma \right) {{\it p_y}
}^{2}}{{x}^{2}+\frac{1}{4}{y}^{2}+\gamma}}\\
F_2&:=&\frac{1}{4}xy^2H+p_y(xp_y-yp_x).
\end{eqnarray*}

If there exists a fourth linear independent integral which is
quadratic in velocities, then, by \cite{konigs},  there exists a
nontrivial integral linear in velocities. Metrics admitting
integrals linear in velocities are described in
Theorem~\ref{linear}. It is easy to see that the lever curves of
the coefficient $f$ from Theorem~\ref{linear}  are generalized
circles ( i.e. circles or straight lines). Since linear
transformations of $\mathbb{C}$
take generalized circles to generalized circles,
and since not all level curves of the function
${x}^{2}+\frac{1}{4}{y}^{2}+\gamma$ are  generalized circles, the   geodesic flow of
metric from Example 2 admits no integral which is linear in
velocities.  Thus, the dimension of  ${\cal
B}(\mathbb{R}^2,g)$ is precisely three.

It is easy to  see that the  linear combination $A H+ B F_1+ C
F_2$ of the integrals  is nonnegative on the whole
$T^*\mathbb{R}^2$, if and only if $B=C=0$. Indeed, the coefficient
by $p_x^2$ is $\frac{4A+By^2+Cy^2x}{x^2+\frac{y^2}{4}+\gamma}$. It
is nonnegative  at every point  of ${\mathbb{R}^2}$ if and only if
$C=0$, $A\ge 0$, $B\ge 0$.  The coefficient by $p_y^2$ in the
linear combination $A H+ B F_1$ is
$\frac{A-B(x^2+\gamma)}{x^2+\frac{y^2}{4}+\gamma}$. It is
nonnegative  at every point  of $\mathbb{R}^2$ if and only if  $A\ge
0$, $B\le 0$. Thus, the  linear combination $A H+ B F_1+ C F_2$ of
the integrals  is nonnegative on the whole $T^*\mathbb{R}^2$, if and
only if $B=C=0$.

 Since every
metric projectively equivalent to $g$ and not proportional to $g$
gives by Corollary~\ref{15} a nonnegative quadratic integral
independent of the energy integral, there is no metric
projectively equivalent to $g$ and nonproportional to $g$.
Theorem~\ref{examples} is proved.

\end{document}